\documentclass[11pt,leqno]{article}

\usepackage{amsmath,amsfonts,amscd,amssymb,theorem}

\long\def\comment#1\endcomment{}

\comment
\endcomment

\makeatletter
\begingroup
\gdef\th@dotted{\normalfont\itshape
  \def\@begintheorem##1##2{%
        \item[\hskip\labelsep \theorem@headerfont ##1\ ##2.]}%
\def\@opargbegintheorem##1##2##3{%
   \item[\hskip\labelsep \theorem@headerfont ##1\ ##2\ (##3).]}}
\endgroup
\makeatother

\theoremstyle{dotted}

\newtheorem{theorem}{Theorem}[section]
\newtheorem{lemma}[theorem]{Lemma}

\newtheorem{prop}[theorem]{Proposition}
\newtheorem{corr}[theorem]{Corollary}

\makeatletter
\begingroup
\gdef\th@upshape{\normalfont
  \def\@begintheorem##1##2{%
        \item[\hskip\labelsep \theorem@headerfont ##1\ ##2.]}%
\def\@opargbegintheorem##1##2##3{%
   \item[\hskip\labelsep \theorem@headerfont ##1\ ##2\ (##3).]}}
\endgroup
\makeatother

\theoremstyle{upshape}

\newtheorem{defn}[theorem]{Definition}
\newtheorem{remark}[theorem]{Remark}

\makeatletter
\renewcommand{\subsection}{\@startsection{subsection}{2}{0pt}{-3ex
plus -1ex minus -0.2ex}{-2mm plus -0pt minus
-2pt}{\normalfont\bfseries}} \makeatother

\makeatletter
\@addtoreset{equation}{section}
\makeatother

\newcommand{\cntrct}                
{\hspace{2pt}\raisebox{1pt}{\text{$\lrcorner$}}\hspace{2pt}}

\newcommand{\proof}[1][Proof.]{\smallskip\noindent{\em #1}}
\def\endproof{\hfill\ensuremath{\square}\par\medskip}

\def\eqref#1{\thetag{\ref{#1}}}

\let\latexref=\ref
\def\ref#1{{\normalfont{\latexref{#1}}}}

\newcommand{\wt}{\widetilde}
\newcommand{\wh}{\widehat}
\newcommand{\whtimes}{\widehat{\times}}
\newcommand{\whotimes}{\operatorname{\wh{\otimes}}}

\newcommand{\ratto}{\dasharrow}        

\newcommand{\idot}{{\:\raisebox{1pt}{\text{\circle*{1.5}}}}}
\newcommand{\hdot}{{\:\raisebox{3pt}{\text{\circle*{1.5}}}}}
\newcommand{\lotimes}{\overset{{\sf L}}{\otimes}}
\newcommand{\compl}{\wt}
\newcommand{\compa}{\overline}

\newcommand{\tw}{ {(1)} }
\newcommand{\Fr}{{\sf Fr}}
\newcommand{\ad}{\operatorname{\sf ad}}
\newcommand{\cchar}{\operatorname{\sf char}}
\newcommand{\Pic}{\operatorname{\sf Pic}}

\newcommand{\Spec}{\operatorname{Spec}}
\newcommand{\Proj}{\operatorname{Proj}}
\newcommand{\sspec}{\operatorname{{\cal S}{\it pec}}}

\newcommand{\Shv}{\operatorname{Shv}}
\newcommand{\fmod}{{\amod^{\text{{\tt\tiny fg}}}}}
\newcommand{\amod}{{\text{\rm -mod}}}

\newcommand{\C}{{\mathbb C}}

\newcommand{\Q}{{\mathbb Q}}
\newcommand{\Z}{{\mathbb Z}}

\newcommand{\aA}{{\mathbb A}}

\newcommand{\Hh}{{\mathbb H}}

\newcommand{\Pp}{{\mathbb P}}
\newcommand{\PP}{{\bf P}}

\newcommand{\X}{{\mathfrak X}}

\newcommand{\A}{{\cal A}}

\newcommand{\D}{{\cal D}}

\newcommand{\F}{{\cal F}}
\newcommand{\M}{{\cal M}}

\newcommand{\LL}{{\cal L}}
\newcommand{\T}{{\cal T}}

\newcommand{\Y}{{\cal Y}}

\newcommand{\E}{{\cal E}}
\newcommand{\K}{{\cal K}}

\newcommand{\G}{{\mathcal G}}

\newcommand{\calo}{{\cal O}}

\newcommand{\m}{{\mathfrak m}}

\newcommand{\hh}{{\cal H}}
\newcommand{\gm}{{\mathbb{G}_m}}

\renewcommand{\phi}{\varphi}

\renewcommand{\dim}{\operatorname{\sf dim}}
\newcommand{\codim}{\operatorname{\sf codim}}

\newcommand{\id}{\operatorname{\sf id}}
\newcommand{\rk}{\operatorname{\sf rk}}
\newcommand{\Id}{\operatorname{{\sf Id}}}

\newcommand{\Hom}{\operatorname{Hom}}
\newcommand{\Ext}{\operatorname{Ext}}

\newcommand{\End}{\operatorname{End}}
\newcommand{\eend}{\operatorname{{\cal E}{\it nd}}}
\newcommand{\Aut}{\operatorname{{\sf Aut}}}

\newcommand{\RHom}{\operatorname{RHom}}
\newcommand{\Rhom}{\operatorname{{\bf R}{\cal H}{\it om}}}

\newcommand{\Ker}{\operatorname{{\sf Ker}}}
\renewcommand{\Im}{\operatorname{{\sf Im}}}

\newcommand{\Supp}{\operatorname{\sf Supp}}

\title{Derived equivalences by quantization}

\author{D. Kaledin\thanks{Partially supported by CRDF grants
RM1-2354-MO02 and RUM1-2694.}}

\date{{\em To Joseph Bernstein, with admiration, on the occasion of
his 60th birthday}}

\begin{document}

\maketitle

\tableofcontents

\section*{Introduction}

The goal of this paper is to find generalization of the so-called
{\em McKay equivalence} of derived categories, as described in
M. Reid's well-known preprint \cite{R}. We briefly recall the
setup. One considers a vector space $V$ over a field of
characteristic $0$ and a finite subgroup $G \subset
SL(V)$. Moreover, one assumes given a smooth crepant resolution $X$
of the quotient variety $Y=V/G$ (crepant in this context means that
the canonical bundle $K_X$ is trivial). In these assumptions, the
{\em McKay correspondence} predicts certain numerical data of $X$,
such as its Betti numbers, starting from the combinatorics of
$G$-action on $V$. This was described in \cite{R}, with precise
conjectures (which were later proved, \cite{ba},
\cite{dl}). However, M. Reid went further: in trying to explain the
geometry underlying his numerical predictions, he proposed a series
of statements, each one stronger than the preceding one. The
strongest of them was the following: one conjectures that the
derived category $\D^b_{coh}(X)$ of coherent sheaves on $X$ is
equivalent to the derived category of $G$-equivariant coherent
sheaves on $V$.

Five years ago, a spectacular proof of this conjecture was given in
\cite{bkr}, under assumption $\dim V = 3$, and for some specific
crepant resolution $X$ (whose existence the authors also
prove). Since then, there has been a lot of progress in $\dim 3$,
and some partial results on adapting the methods of \cite{bkr} to
some cases of higher dimension (\cite{go}). However, to the best of
our knowledge, the only relatively general result in higher
dimension was obtained very recently by the author jointly with
R. Bezrukavnikov, \cite{BK2}. We establish the McKay equivalence in
arbitrary dimension, but under one additional assumption: we require
$V$ to be a {\em symplectic} vector space, and we want the finite
group $G$ to preserve the symplectic form.

But there was a different development still in dimension $3$. It was
realized by T. Bridgeland \cite{brid} that the methods of \cite{bkr}
can work in a more general situation. One still considers a crepant
resolution $X$ of a singular affine algebraic variety $Y$, but one
no longer requires $Y$ to be a quotient singularity. In \cite{brid},
Bridgeland considers instead a so-called {\em small contraction} $X
\to Y$ of a smooth $3$-dimensional manifold $X$ with trivial $K_X$
-- that is, he assumes given a proper birational map $X \to Y$ whose
only exceptional fibers are curves. Bridgeland's results were
extended and re-cast in a different language by M. Van den Bergh
\cite{vdb1}, and in this form, they are very similar to the McKay
equivalence: the derived category $D^b_{coh}(X)$ of coherent sheaves
on $X$ is shown to be equivalent to the derived category of
finite-generated left modules over a certain non-commutative algebra
$R$. The algebra $R$ has a structure formalized by Van den Bergh
\cite{vdb2} under the name of a {\em non-commutative resolution} of
the affine variety $Y$; in particular, $R$ has a large center, which
is identified with the algebra $A=H^0(Y,\calo_Y)$ of functions on
$Y$, and generically on $Y$ -- that is, after tensoring with the
fraction field of $A$ -- the algebra $R$ is a matrix algebra. It
also enjoys some other nice properties, such as finite global
homological dimension.

We note that the McKay equivalence can be stated in exactly the same
way, with $R$ being the so-called {\em smash product algebra}
$S^\hdot(V^*) \# G$ (see e.g. \cite{BK2}); the category of finitely
generated left $R$-modules is then immediately seen to be equivalent
with the category of $G$-equivariant coherent sheaves on $V$. All in
all, it seems that Van den Bergh's non-commutative resolution
picture is the proper framework for generalizing McKay equivalence,
and it is this picture that one should try to find in higher
dimensions.

\medskip

In this paper, we do this under the same additional assumption as in
\cite{BK2}: we assume given a {\em symplectic} resolution $X$ of a
normal irreducible affine variety $Y$, and we construct a
non-commutative resolution $R$ of the variety $Y$ and an equivalence
between $D^b_{coh}(X)$ and the bounded derived category of finitely
generated left $R$-modules. Unfortunately, we can only do this
locally on $Y$ -- that is, we fix a point $y \in Y$, and in the
course of our construction we may have replace $Y$ with an \'etale
neighborhood of the point $y$. However, we impose no additional
restrictions on $X$. We also prove that if $X$, $X'$ are two
different symplectic resolutions of the same variety $Y$, then,
locally on $Y$, the derived categories $D^b_{coh}(X)$ and
$D^b_{coh}(X')$ are equivalent (this generalizes the particular case
of \cite[Section 3, Conjecture]{BO1} proved by Y. Kawamata in
\cite{Ka}: in Kawamata's language, ``$K$-equivalence implies
$D$-equivalence''). As in \cite{BK2}, our main technical tool is
reduction to positive characteristic and applying the Fedosov
quantization procedure, which has recently been worked out in
positive characteristic in \cite{BK3}.

\medskip

The paper is organized as follows. In Section~\ref{intr} we give the
precise statements of our results, and also do some reformulation;
in particular, we prove that the main Theorem is equivalent to $X$
having a so-called {\em tilting generator} $\E$. In
Section~\ref{prem}, we recall additional material needed for the
construction: in Subsection~\ref{tw}, we introduce a certain
one-parameter deformation of the symplectic manifold $X$ called the
{\em twistor deformation}, and in Subsection~\ref{qua}, we recall
the necessary facts about reduction to positive characteristic and
quantization over positive characteristic fields. At this point, we
are able to construct a tilting vector bundle $\E$ on $X$. To prove
the main Theorem, it remains to show that $\E$ is a generator; this
is done in Section~\ref{d.aff}, after preliminary technical
estimates are established in Section~\ref{est}. Finally, in
Section~\ref{app} we give some applications and generalizations of
our main result.

\medskip

One final remark is perhaps in order. In \cite{BK2}, when dealing
with symplectic resolutions of quotient singularities, not only did
we establish an equivalence $D^b_{coh}(X) \cong D^b(R\fmod)$, but we
were actually able to compute the algebra $R$. In this paper, we do
not try to do this. However, we expect that this is a meaningful
thing to do, and that the resulting algebras should be related to
the theory of quantum groups. For more details, see Remark~\ref{qg}.

\subsection*{Acknowledgements.} This paper grew out of a long joint
project with R. Bezrukavnikov, which also produced \cite{BK},
\cite{BK2} and \cite{BK3}; many of the ideas in the present paper
are actually his ideas. I have also benefitted greatly from
discussions with A. Braverman and D. Kazhdan. In particular,
A. Braverman contributed a lot to the material in
Subsection~\ref{pure.sub}. Many discussions with V. Fock,
V. Ginzburg, D. Huybrechts, K. Kremnitzer, A. Kuznetsov, M. Lehn and
N. Markarian were very helpful. A. Bondal, D. Matsushita,
Y. Namikawa and D. Orlov have read a first version of the paper and
made several important suggestions. A. Kuznetsov read the first
version quite thoroughly and found a sizable number of bugs, for
which I am extremely grateful. Part of this research was carried out
at the Hebrew University of Jerusalem in Israel; it is a pleasure to
thank this wonderful institution for an opportunity to visit and to
do some work.

\section{Statements.}\label{intr}

Let $X$ be a smooth manifold -- that is, a regular finite-type
scheme -- over some field $k$. We will say that $X$ is {\em convex}
if it is equipped with a projective birational map $\pi:X \to Y$
onto a normal irreducible affine algebraic variety $Y$ of finite
type over $k$ (we note that under these assumptions we must have
$Y=\Spec H^0(X,\calo_X)$, so that the existence of $Y$ and $\pi$ is
a condition on $X$, not some new data). We start with the following
general definition.

\begin{defn}\label{tilt}
A coherent sheaf $\E$ on $X$ is called a {\em tilting generator} of
the bounded derived category $D^b_{coh}(X)$ of coherent sheaves on
$X$ if the following holds:
\begin{enumerate}
\item The sheaf $\E$ is a tilting object in $D^b_{coh}(X)$ -- that
  is, for any $i \geq 1$ we have $\Ext^i(\E,\E)=0$
\item The sheaf $\E$ generates the derived category $D^-_{coh}(X)$
  of complexes bounded from above -- that is, if for some object $\F
  \in D^-_{coh}(X)$ we have $\RHom^\hdot(\E,\F)=0$, then $\F=0$.
\end{enumerate}
\end{defn}

\begin{lemma}\label{tilt.lemma}
Assume that $X$ is convex, let $\E$ be a tilting generator of the
derived category $D^b_{coh}(X)$, and denote $R=\End(\E)$. Then the
algebra $R$ is left-Noetherian, and the correspondence $\F \mapsto
\RHom^\hdot(\E,\F)$ extends to an equivalence
\begin{equation}\label{equi}
D^b_{coh}(X) \to D^b(R\fmod)
\end{equation}
between the bounded derived category $D^b_{coh}(X)$ of coherent
sheaves on $X$ and the bounded derived category $D^b(R\fmod)$ of
finitely generated left $R$-modules.
\end{lemma}

\proof{} Since $X$ is convex, the algebra $R$ is a module of finite
type over the commutative algebra $H^0(X,\calo_X)$ of global
functions on $X$, and this commutative algebra by definition is the
algebra of functions on an affine variety $Y$ of finite type over
$k$. Hence $H^0(X,\calo_X)$ is Noetherian, and $R$ is left (and
right, and two-sided) Noetherian. In particular, the category
$R\fmod$ is abelian, and $D^b(R\fmod)$ is well-defined.

The functor $a:\F \mapsto \Hom(\E,\F)$ is a left-exact functor from
the abelian category of coherent sheaves on $X$ to the abelian
category of finitely generated $R$-modules. Its derived functor is
therefore a well-defined functor $A:D^b_{coh}(X) \to
D^b(R\fmod)$. It is easy to see that $a$ has a right-exact
left-adjoint functor $b:M \mapsto M \otimes_R \E$; its derived
functor $B:D^-(R\fmod) \to D^-_{coh}(X)$ is adjoint to $A$. The
composition
$$
A \circ B:D^-(R\fmod) \to D^-(R\fmod)
$$
sends $R$ to $\Rhom^\hdot(\E,\E)$; since $\E$ is tilting, we have $A
\circ B(R)=R$. But $R$ generates $D^-(R\fmod)$. Therefore
$B:D^-(R\fmod) \to D^-_{coh}(X)$ is fully faithful. Moreover, for
any $\F \in D^-_{coh}(X)$, the cone of the adjunction map $B(A(\F))
\to \F$ is annihilated by $A$, which by
Definition~\ref{tilt}~\thetag{ii} means that $\F \cong
B(A(\F))$. Therefore $A$ and $B$ are mutually inverse equivalences
between $D^-_{coh}(X)$ and $D^-(R\fmod)$. This in particular means
that $B$ is adjoint to $A$ {\em both on the right and on the left};
since $A$ sends $D^b_{coh}(X)$ into $D^b(R\fmod)$ and has bounded
cohomological dimension, this implies that $B$ sends $D^b(R\fmod)$
into $D^b_{coh}(X)$, so that $A$ and $B$ also induce mutually
inverse equivalences between $D^b_{coh}(X)$ and $D^b(R\fmod)$.
\endproof

\begin{remark}
In Definition~\ref{tilt}, we require $\E$ to be a coherent sheaf on
$X$, not just an arbitrary object in $D^b_{coh}(X)$ or
$D^-_{coh}(X)$. This might not be strictly necessary for
Lemma~\ref{tilt.lemma}, but it simplifies the proof, and this level
of generality is sufficient for our purposes. In our applications,
$\E$ will in fact be not just a sheaf, but a vector bundle.
\end{remark}

Assume now that the base field $k$ has characteristic $0$, and that
$X$ is symplectic -- that is, we are given a non-degenerate closed
$2$-form $\Omega \in H^0(X,\Omega^2_X)$. Moreover, assume that $X$
is convex, with $Y = \Spec H^0(X,\calo_X)$. The main result of the
paper is the following.

\begin{theorem}\label{main}
Under the assumptions above, for any point $y \in Y$ there exists an
\'etale neighborhood $U_y \to Y$ such that the pullback $X_y=X
\times_Y U_y$ admits a tilting generator $\E$. Moreover, this
tilting generator is in fact a vector bundle on $X_y$.
\end{theorem}

\begin{remark}
The standard $t$-structure on $D^b(R\fmod)$ gives by \eqref{equi} a
non-standard $t$-structure on $D^b_{coh}(X)$ and defines a notion of
a perverse coherent sheaf on $X$ (it is this non-standard
$t$-structure, rather than the equivalence \eqref{equi}, that was
discovered in \cite{brid}). Our claim that the tilting generator
$\E$ is a vector bundle is then obviously equivalent to the fact
that every skyscraper coherent sheaf on $X$ is perverse.
\end{remark}

We also prove the following.

\begin{theorem}\label{K=>D}
Assume given two smooth projective resolutions $X$, $X'$ of a normal
affine irreducible variety $Y$, assume that the canonical bundles
$K_X$, $K_{X'}$ are trivial, and assume moreover that $X$ admits a
closed non-degenerate $2$-form $\Omega \in H^0(X,\Omega^2_X)$. Then
every point $y \in Y$ admits an \'etale neighborhood $U_y \to Y$
such that the derived categories $D^b_{coh}(X \times_Y U_y)$ and
$D^b_{coh}(X' \times_Y U_y)$ are equivalent.
\end{theorem}

This is a very particular case of \cite[Conjecture 1.2]{Ka}, which
in turn essentially goes back to \cite[Section 3]{BO1} (see also
\cite{BO2}). In the language of Y. Kawamata, ``$K$-equivalence
implies $D$-equivalence'' (for symplectic resolutions, and locally
over the base). In the original, slightly less precise language of
A. Bondal and D. Orlov, ``if two smooth symplectic varieties are
related by a flop, their derived categories are equivalent''.

The need to pass to an \'etale neighborhood of a point is
unfortunate, but, seemingly, unavoidable in the general situation.
In practice, the problem can sometimes be alleviated by presence of
an additional structure. We prove one result of this sort.

\begin{defn}
An action of the multiplicative group $\gm$ on an affine scheme $Y$
with fixed closed point $y \in Y$ is said to have {\em positive
weights} if the weights of the $\gm$-action are non-negative on the
function algebra $H^0(Y,\calo_Y)$, and strictly positive on the
maximal ideal $\m \subset Y$ which defines the point $y$.
\end{defn}

\begin{theorem}\label{gm}
\begin{enumerate}
\item Assume that a smooth symplectic scheme $X$ is projective over
an affine variety $Y$, and assume that $Y$ admits an action of the
group $\gm$ with fixed closed point $y \in Y$ and positive
weights. Then the $\gm$-action extends to a $\gm$-action on $X$.
\item Assume that a scheme $X$ is projective over an affine variety
  $Y$, and assume that $Y$ admits an action of the group $\gm$ with
  fixed closed point $y \in Y$ and positive weights which lifts to a
  $\gm$-action on $X$. Assume also that for some \'etale neighborhood $U_y$
  of the point $y \in Y$, the pullback $X \times_Y U_y$ admits a
  tilting generator $\E_y$. Then $\E_y$ is obtained by pullback from
  a $\gm$-equivariant tilting generator $\E$ on $X$.
\end{enumerate}
\end{theorem}

Finally, we prove one cohomological consequence of the existence of
tilting generators.

\begin{theorem}\label{diag}
Assume that a smooth manifold $X$ is projective over an affine local
Henselian scheme $Y/k$ and admits a tilting generator $\E$. Then the
structure sheaf $\calo_\Delta$ of the diagonal $\Delta \subset X
\times X$ admits a finite resolution by vector bundles of the form
$\E_i \boxtimes \F_i$, where $\E_i$, $\F_i$ are some vector bundles
on $X$.
\end{theorem}

\begin{corr}\label{pure}
Assume that a smooth manifold $X$ is projective over an affine
scheme $Y$, and let $F \subset X$ be the fiber over a closed point
$y \in Y$. Assume that $Y$ admits a positive-weight $\gm$-action
that fixes $y \in Y$, and assume that $X$ admits a tilting generator
$\E$. Then the cohomology groups $H^\hdot(F)$ of the scheme $F$ are
generated by classes of algebraic cycles.
\end{corr}

In this Corollary we are deliberately vague as to what particular
cohomology groups $H^\hdot(F)$ one may take. In fact, every
cohomology theory with the standard weight formalism will suffice;
in particular, the statement is true for $l$-adic cohomology and for
analytic cohomology when the base field $k$ is $\C$.

\section{Preliminaries.}\label{prem}

\subsection{Twistor deformations.}\label{tw}
Assume given a smooth manifold $X$ over a field $k$ equipped with a
non-degenerate closed $2$-form $\Omega \in H^0(X,\Omega^2_X)$ (from
now on, we will call such a form a symplectic form). Assume also
given a projective birational map $\pi:X \to Y$ onto a normal
irreducible affine scheme $Y$ of finite type over $k$. Choose a line
bundle $L$ on $X$. Denote $S=\Spec k[[t]]$, the formal disc over
$k$, and let $o \in S$ be the special point (given by the maximal
ideal $tk[[t]] \subset k[[t]]$).

\begin{defn}
By a {\em twistor deformation} $Z$ associated to the pair $\langle
X,L \rangle$ we will understand a triple of a smooth scheme $\X/S$,
a line bundle $\LL$ on $\X$ and a symplectic form $\Omega_Z$ on the
total space $Z$ of the $\gm$-torsor associated to $\LL$ (in other
words, on the total space of $\LL$ without the zero section) such
that
\begin{enumerate}
\item The form $\Omega_Z$ is $\gm$-invariant, and the map $\rho:Z
  \to S$ is the moment map for the $\gm$-action on $Z$ -- that is,
  we have $\Omega_Z \cntrct \xi_0 = \rho^*dt$, where $\xi_0 \in
  H^0(Z,T_Z)$ is the infinitesimal generator of the $\gm$-action.
\item The restriction $\langle \X_o, \LL_o \rangle$ of the pair
  $\langle \X,\LL \rangle$ to the special point $o \in S$ is
  identified with the pair $\langle X, L \rangle$, and the
  restriction of the form $\Omega_Z$ to the special fiber $Z_o
  \subset Z$ coincides under this identification with the pullback
  of the given form $\Omega \in H^0(X,\Omega^2_X)$.
\end{enumerate}
\end{defn}

Assume that the base field $k$ has characteristic $0$. Then for any
$i \geq 1$ we have $H^i(X,K_X)=0$ by the Grauert-Riemenschneider
Vanishing Theorem, and since the top power of the symplectic form
$\Omega$ trivializes the canonical bundle $K_X$, this implies that
$H^i(X,\calo_X)=0$ for $i \geq 1$. Therefore the variety $X$ falls
within the assumptions of the paper \cite{K1}.

\begin{lemma}\label{tw.lemma}
For any line bundle $L$ on $X$, there exists a twistor deformation
$\langle \X,\LL,\Omega_Z\rangle$ associated to the pair $\langle X,L
\rangle$. Moreover, $\X$ is projective over $\Y = \Spec
H^0(\X,\calo_\X)$, while $\Y$ is normal and flat over $S$.
\end{lemma}

\proof{} This is a particular case of \cite[Theorem 1.4]{K1}; the
second statement is \cite[Theorem 1.5]{K1}.
\endproof

There is also a certain uniqueness statement in \cite[Theorem
1.4]{K1}, but we will not need it. What is important is that the
construction is sufficiently functorial. This allows to prove the
following.

\begin{defn}\label{tw.ex}
A twistor deformation $\langle\X,\LL,\Omega_Z\rangle$ is called {\em
exact} if the symplectic form $\Omega_Z$ is exact,
$\Omega_Z=d\alpha_Z$, and moreover the $1$-form $\alpha$ is
$\gm$-equivariant and satisfies $\alpha \cntrct \xi_0=\rho^*t$,
where $\xi_0$ is the infinitesimal generator of the $\gm$-action on
$Z$, $t$ is the coordinate on $S = \Spec k[[t]]$, and $\rho:Z \to S$
is the natural projection.
\end{defn}

\begin{lemma}
Assume that the symplectic form $\Omega$ on $X$ is exact. Then any
twistor deformation $\langle\X,\LL,\Omega_Z\rangle$ associated to
the pair $\langle X,L \rangle$ by Lemma~\ref{tw.lemma} is also
exact.
\end{lemma}

\proof{} Since $\Omega$ is non-degenerate, $\alpha = \Omega \cntrct
\xi$ for some vector field $\xi \in H^0(X,T_X)$ on $X$. By the
Cartan Homotopy formula, $d\alpha = \Omega$ is equivalent to
$\LL_\xi(\Omega)=\Omega$, where $\LL_\xi$ is the Lie derivative
along $\xi$. Since $H^1(X,\calo_X)=0$, the line bundle $L$ can be
made equivariant with respect to the vector field $\xi$ -- indeed,
the Atiyah extension
$$
\begin{CD}
0 @>>> \calo_X @>>> \E @>>> \T_X @>>> 0
\end{CD}
$$
associated to $L$ must split after restricting to the section
$\xi:\calo_X \to \T_X$ of the tangent bundle $\T_X$. Thus $\xi$ is
an infinitesimal automorphism of the pair $\langle X,L \rangle$, and
it dilates the symplectic form. In terms of the Poisson structure on
$X$, this is equivalent to \cite[(1.2)]{K1}.  By \cite[Proposition
1.5]{K1} the vector field $\xi$ extends to a $\gm$-invariant vector
field $\xi_Z$ on $Z$ such that $\xi_Z(\rho^*t)=-t$. By the Cartan
Homotopy Formula, $\Omega_Z = d(\Omega_Z \cntrct \xi_Z)$, so that we
can take $\alpha_Z=\Omega_Z \cntrct \xi_Z$. Since both $\xi_Z$ and
$\Omega_Z$ are $\gm$-invariant, so is $\alpha_Z$. Moreover, we have
$$
\alpha_Z \cntrct \xi_0 = \Omega_Z \cntrct (\xi_Z \wedge \xi_0) = -
\Omega_Z \cntrct (\xi_0 \wedge \xi_Z) = -\rho^*dt \cntrct \xi_Z =
-\xi_Z(\rho^*t) = \rho^*t,
$$
which proves that $Z$ is indeed exact in the sense of
Definition~\ref{tw.ex}.
\endproof

We note that by \cite[Corollary 2.8]{K2}, the symplectic form
$\Omega$ is always exact over a formal neighborhood of any closed
point $y \subset Y$. Finally, we will need the following.

\begin{lemma}\label{1-1}
Assume that the line bundle $L$ on $X$ is ample, and consider the
twistor deformation $\langle \X,\LL,\Omega_Z\rangle$ associated to
the pair $\langle X,L \rangle$ by Lemma~\ref{tw.lemma}. Let $\wh{A}
= H^0(\X,\calo_{\X})$, $\Y = \Spec \wh{A}$, and let $\pi:\X \to \Y$
be the natural map. Then the map $\pi$ is projective, and one-to-one
over the complement $S \setminus \{o\}$. Moreover, if $Y$ is the
spectrum of a Henselian local $k$-algebra, so that $\wh{A}$ is a
local $k$-algebra with maximal ideal $\m \subset \wh{A}$, then there
exists a finitely generated $k$-subalgebra $\wt{A} \subset \wh{A}$
such that
\begin{enumerate}
\item the $t$-adic completion of the Henselization of the algebra
  $\wt{A}$ in $\m \cap \wt{A} \subset \wt{A}$ coincides with
  $\wh{A}$, and
\item all the data $\langle \X,\LL,\Omega_Z \rangle$ are defined
over $\wt{A}$.
\end{enumerate}
\end{lemma}

\proof{} (The argument uses a standard trick which probably goes
back to \cite{F}. Compare \cite[Propositon 4.1]{H1} and \cite[Claim
3 in the proof of Theorem 2.2]{N1}.)

The generic fiber $\Y_\eta = \Y \times_S \eta$ is open in the normal
scheme $\Y$ (indeed, it is the complement to the closed special
fiber $Y \subset \Y$). Since the scheme $\Y$ is normal, the open
subscheme $\Y_\eta \subset \Y$ is also normal.  Since the map
$\pi:\X \to \Y$ is birational, it suffices to prove that it is
finite over the generic point $\eta \in S$. Moreover, by
construction the map $\pi:\X_\eta \to \Y_\eta$ is
projective. Therefore by \cite[IV, Th\'eor\`eme 8.11.1]{EGA} it
suffices to prove that it is quasifinite -- in other words, that its
fibers do not contain any proper curves.

Let $\iota:C \to \X_\eta$ be an arbitrary map from a proper curve
$C/\eta$ to $\X_\eta$. Replacing $C$ with its normalization, we can
assume that the curve $C/\eta$ is connected and smooth. The
$\gm$-equivariant symplectic form $\Omega_Z$ induces by descent a
relative symplectic form $\Omega_\X \in H^0(\X,\Omega^2_{\X/S})$,
and by \cite[Lemma 1.5]{K1}, its de Rham cohomology class
$[\Omega_\X]$ is equal to $[\Omega]+c_1(L)t$, where $c_1(L)$ is the
first Chern class of the line bundle $L$ (here we have identified
$H^2_{DR}(\X/S)$ with $H^2_{DR}(X)[[t]]$ by means of the Gauss-Manin
connection). Since $C$ is a curve, $\iota^*\Omega_\X$ is trivial;
therefore
$$
0= \iota^*[\Omega_\X] = \iota^*[\Omega]+t\iota^*c_1(L).
$$
Differentiating this equality with respect to $t$, we obtain
$\iota^*c_1(L) =0$. Since $L$ is an ample line bundle, its extension
$\LL$ is also ample. The first Chern class $c_1(\LL)$ is constant
with respect to the Gauss-Manin connection; therefore
$\iota^*c_1(\LL) = \iota^*c_1(L)=0$, and since $\LL$ is ample, this
is possible only if $\iota:C \to \X_\eta$ maps the curve $C$ to a
point.

Finally, to show that the formal scheme $\Y$ comes from an algebraic
variety $\Spec \wt{A}$, we note that since $\X$ is smooth and
$\pi:\X \to \Y$ is bijective outside of the special fiber $Y \subset
\Y$, the scheme $\Y$ is smooth outside of $Y \subset \Y$. Then the
existence of a finitely generated subalgebra $\wt{A} \subset \wh{A}$
which satisfies \thetag{i} is insured by \cite[Theorem
3.9]{A}. Localizing $\wt{A}$ if necessary, we can also achieve
\thetag{ii}.
\endproof

\subsection{Quantization.}\label{qua}
We will now assume that the base field $k$ is perfect and has odd
positive characteristic $p$. In this section, we prove the main
quantization result which we need in the paper,
Proposition~\ref{qua.prop}; the proof depends heavily on notions
introduced in \cite[Section 1]{BK3} (although the statement ought to
be comprehensible to a reader who is not familiar with that paper,
and it can be used as a black box).

For any scheme $Z$ over $k$, we denote by $Z^\tw$ the twist of the
scheme $Z$ with respect to the Frobenius map $k \to k$. As a
topological space, $Z^\tw$ coincides with $Z$, and if the scheme $Z$
is reduced, the structure sheaf $\calo_{Z^\tw}$ is canonically
identified with the subalgebra $\calo_Z^p \subset \calo_Z$ generated
by $p$-th powers of functions on $Z$. On level of schemes, the
embedding $\calo_{Z^\tw} \cong \calo_Z^p \subset \calo_Z$ is the
Frobenius map $\Fr_Z:Z \to Z^\tw$.  In particular, if $S = \Spec
k[[t]]$, then $S^\tw$ is $\Spec k[[t^p]]$, and we have the Frobenius
map $\Fr_S:S \to S^\tw$.

Consider the power series algebra $k[[t,h]]$ in two variables $t$,
$h$, and denote $S_h=\Spec k[[t,h]]$. Define a map $s:k[[t^p]] \to
k[[t,h]]$ by setting $s(t^p) = t^p + th^{p-1}$. By definition, $S =
\Spec k[[t]]$ is canonically embedded into $S_h$; the map $s:S_h \to
S^\tw$ extends the Frobenius map $\Fr_S:S \to S^\tw$ to $S_h \supset
S$. This extension is {\em not} the obvious one, and it behaves
differently. In particular, the fiber $s^{-1}(o) \subset S_h$ over
the special point $o \in S^\tw$ is the union of $p$ formal lines
$t=ah$, $a \in \Z/p\Z \subset k$ in the formal affine plane
$S_h$. The natural projection $\sigma:s^{-1}(o) \to \Spec k[[h]]$ is
the Artin-Schreier covering.

Assume now given a symplectic manifold $X/k$ equipped with a
projective birational map $\pi:X \to Y$ onto a normal irreducible
affine variety $Y/k$; moreover, assume given a line bundle $L$ on
$X$ and a twistor deformation $\langle\X,\LL,\Omega_Z\rangle$ of the
pair $\langle X,L \rangle$. Denote by $\X^\tw_h$ the relative
spectrum
\begin{equation}\label{xh}
\X^\tw_h = \sspec (\X,\calo^p_{\X} \whotimes_{k[[t^p]]} k[[t,h]]),
\end{equation}
where $k[[t^p]]$ is embedded into $k[[t,h]]$ by means of the map
$s$, and $\whotimes$ stands for $h$-adic completion of the tensor
product. The scheme $\X^\tw_h$ is regular and projective over
\begin{equation}\label{yh}
\Y_h^\tw = \Spec H^0\left(\X^\tw_h,\calo^p_{\X^\tw_h}\right) = \Spec
\left(H^0(\X,\calo_\X) \whotimes_{k[[t^p]]} k[[t,h]]\right).
\end{equation}
Denote by $\X^\tw_o \subset \X^\tw_h$ the special fiber of the
natural flat projection $\X^\tw_h \to \Spec k[[h]]$. As a
topological space, $\X^\tw_o$ coincides with $\X$, while the
structure sheaf is given by $\calo_\X^p \otimes_{k[[t^p]]} k[[t]]
\subset \calo_\X$. In particular, $\calo_\X$ is a coherent sheaf on
$\X^\tw_o$, and its restriction $\calo_\X/t$ to $X^\tw \subset
\X^\tw_o$ coincides with the direct image $\Fr_*\calo_X$ under the
Frobenius map (thus $\calo_\X$ is in fact a vector bundle of rank
$p^{\dim X}$ on $\X^\tw_o$).

\begin{prop}\label{qua.prop}
Assume that the twistor deformation $\langle \X,\LL,\Omega_Z
\rangle$ is exact in the sense of Definition~\ref{tw.ex}, and assume
that we have $H^i(\X,\calo_\X)=0$ for $i \geq 1$. Then there exists
a coherent sheaf $\calo_h$ of algebras on $\X^\tw_h$, flat over
$k[[h]]$, such that the restriction $\calo_h|_{\X^\tw_o}$ to
$\X^\tw_o \subset \X^\tw_h$ is identified with $\calo_\X$, while
over the complement $\X^\tw_h \setminus \X$, the algebra $\calo_h$
is isomorphic to the endomorphism algebra of a vector bundle $\E$.
\end{prop}

\proof{} We freely use the notions from \cite[Section 1]{BK3} -- in
particular, that of a restricted Poisson algebra, \cite[Definition
1.8]{BK3}. Recall that a restricted structure on a Poisson algebra
$A$ is defined by a non-additive ``restricted power'' operation $x
\mapsto x^{[p]}$ on $A$ compatible with the multiplication and the
Poisson bracket in a certain specified way; if the Poisson bracket
on $A$ is trivial, then the restricted power must be additive, and
it must be a ``Frobenius-derivation'' -- that is, we have
$(ab)^{[p]} = a^{[p]}b^p + a^p b^{[p]}$.  Let $K_0:k[[t]] \to
k[[t]]$ be the (unique) Frobenius-derivation such that $K_0(t) =
t$. Let $\alpha_Z$ be the $1$-form on $\X$ whose existence is
required by Definition~\ref{tw.ex}. By \cite[Proposition 2.6]{BK3},
setting
$$
f^{[p]}= H_f^{[p]} \cntrct \alpha_Z - H_f^{p-1}(H_f \cntrct
\alpha_Z),
$$
defines a restricted Poisson structure on $Z$ -- here $H_f$ is the
Hamiltonian vector field associated to the function $f$, and
$H_f^{[p]}$ is its restricted $p$-th power with respect to usual
restricted Lie algebra structure on the Lie algebra of vector
fields. Since $\alpha_Z$ is $\gm$-invariant, this restricted Poisson
structure is $\gm$-equivariant; in particular, it descends to a
restricted Poisson structure on the quotient $\X = Z/\gm$. Moreover,
since $H_{\rho^*t} = \xi_0$, the differential of the $\gm$-action,
we have
$$
(\rho^*t)^{[p]}= \xi_0^{[p]} \cntrct \alpha_Z - \xi_0^{p-1}(\xi_0
\cntrct \alpha_Z).
$$
Since $\xi_0$ is the differential of an action of the multiplicative
group, we have $\xi_0^{[p]}=\xi_0$. By assumption $\xi_0 \cntrct
\alpha_Z = \rho^*t$; in particular, it is $\gm$-invariant. We
conclude that $(\rho^*t)^{[p]} = \xi_0 \cntrct \alpha_Z =
\rho^*t$. This means that the restricted Poisson structure on $\X/S$
is compatible with the Frobenius-derivation $K_0$ in the sense of
\cite[Corollary 1.13]{BK3}.

We now notice that $B=k[[t,h]]$ has a natural structure of a
quantization base in the sense of \cite[Definition 1.15]{BK3}. To
define it, one considers the map $s_0:k[[t]] \to k[[h,t]]$,
$s_0(t)=t$, and notices that for any $f \in k[[t]]$, the difference
$s(f^p)-s_0(f)^p \in k[[h,t]]$ is divisible by $h^{p-1}$; thus there
is a unique map $K:k[[t] \to k[[t,h]]$ such that $K(f)h^{p-1} +
s_0(f)^p = s(f^p)$, and we extends it to $k[[t,h]]$ by setting $K(f)
= K(f \mod h)$. We note that $K(t)=t$; therefore $K(f) = K_0(f) \mod
h$ for any $f \in k[[t]]$, and we can apply \cite[Theorem 1.23]{BK3}
to the restricted Poisson structure on $\X/S$. By \cite[Proposition
1.24]{BK3}, the result is a sheaf of algebras $\calo_h$ on
$\X^\tw_h$, flat over $k[[h]]$, whose restriction to $\X^\tw_o$ is
identified with $\calo_\X$ -- or actually, a whole set $Q(B,\X)$ of
isomorphism classes of such sheaves. Moreover, on the complement
$\X^\tw_h \setminus \X^\tw_o$, every sheaf $\calo_h$ is an Azumaya
algebra. Finally, the class of this Azumaya algebra in the Brauer
group is also described in \cite[Proposition 1.24]{BK3} and, as
noted there, one can choose an element in $Q(B,\X)$ so that the
corresponding Azumaya algebra is split, $\calo_h \cong \eend(\E)$ on
$\X^\tw_h \setminus \X^\tw_o$ for some vector bundle $\E$.
\endproof

\begin{lemma}\label{tlt}
In the assumptions of Proposition~\ref{qua.prop}, we have
$\Ext^i(\E,\E)=0$ for $i \geq 1$.
\end{lemma}

\proof{} Since $\eend(\E) \cong \calo_h$ on $\X^\tw_h \setminus
\X^\tw_o$, it suffices to prove that $\calo_h$ has no cohomology on
this complement. By base change, it suffices to prove that it has no
cohomology on the whole $\X^\tw_h$. To compute
$H^\hdot(\X^\tw_h,\calo_h)$ one can use the spectral sequence
associated to the $h$-adic filtration on the $k[[h]]$-flat sheaf
$\calo_h$, and since the quotient $\calo_h/h$ is supported on
$\X^\tw_o$, the term $E^1$ of this sequence is
$$
H^i(\X^\tw_o,\calo_h/h)[[h]].
$$
Since $\X^\tw_o = \X$ as a topological space and $\calo_h/h \cong
\calo_\X$, this vanishes for $i \geq 1$ by assumption.
\endproof

\begin{remark}\label{qg}
As we have noted in the Introduction, we make no attempt to describe
the algebra $R=\End(\E)$ obtained by this construction, nor its
lifting to characteristic $0$ which we will construct later in
Section~\ref{d.aff}. However, if one replaces $X$ with a semisimple
algebraic group $G$ considered as a Poisson-Lie group, then the
resulting algebra $R$ is very similar to the dual to the quantum
envelopping algebra specialized at a $p$-th root of unity introduced
by G. Lusztig. The similarity becomes even more pronounced if one
restricts the quantum group to a symplectic leaf in $G$, as done,
for instance, by C. De Concini and C. Procesi in \cite{CP}. It must
be noted, however, that Lusztig -- and consequently, De
Concini-Procesi -- work by explicit computation. Our approach only
works for symplectic $X$, but it makes no use of any group
structure, and it somewhat clarifies the geometric picture. In a
nutshell, what happens is this: after one reduces a Poisson variety
$X$ to positive characteristic, it acquires the Frobenius map $F:X
\to X$. This map of course cannot in general be lifted back to
characteristic $0$. What we prove for symplectic $X$ is that it does
admit a lifting {\em if we first deform the function algebra of $X$
to a quantized function algebra}. The result is a quantum analog of
Frobenius map in characteristic $0$ -- when $X=G$ is a semisimple
group, this is Lusztig's quantum Frobenius map. Motivated by this
analogy, we expect that with appropriate modifications, out
construction should work for at least some Poisson manifolds which
are not symplectic. A natural question is to find the proper
conditions on the Poisson manifold $X$ which allow to construct a
quantum Frobenius map.
\end{remark}

\section{Estimates.}\label{est}

\subsection{Generalities on algebra sheaves.}

We start with some generalities on algebra sheaves. Consider a
scheme $X$ flat over some affine scheme $S=\Spec B$, and assume that
$X$ is equipped with a coherent flat algebra sheaf $\A$. Denote by
$A = H^0(X,\A)$ the algebra of global sections of the sheaf $\A$.
Taking global sections is then a left-exact functor from the
category $\Shv(X,\A)$ of quasicoherent sheaves of left $\A$-modules
on $X$ to the category $A\amod$ of left $A$-modules. We denote this
functor -- for reasons of convenience -- by $\pi^X_*:\Shv(X,\A) \to
A\amod $. It has a left-adjoint right-exact functor
$\pi_X^*:A\amod\to\Shv(X,\A)$ given by $\pi^*(M) = M \otimes_A \A$.

Denote by $R^\hdot\pi^X_*$ and $L^\hdot\pi_X^*$ the derived functors
of $\pi^X_*$ and $\pi_X^*$. Assume that the category $\Shv(X,\A)$
has finite homological dimension; then $R^\hdot\pi^X_*$ and
$L^\hdot\pi_X^*$ are a pair of adjoint functors between the derived
categories $\D^-(X,\A) = \D^-\Shv(X,\A)$ and $\D^-(A) =
\D^-A\amod$. The composition
$$
R^\hdot\pi^X_*L^\hdot\pi_X^*:D^-(A) \to D^-(A)
$$ 
comes equipped with the adjunction map $\Id \to
R^\hdot\pi^X_*L^\hdot\pi_X^*$. It is easy to see that this map is an
isomorphism if and only if $H^i(X,\A)=0$ for $i \geq 1$: indeed,
since every module has a free resolution, it suffices to check that
$A \to R^\hdot\pi^X_*L^\hdot\pi_X^*A$ is an isomorphism, and the
right-hand side by definition coincides with $H^\hdot(X,\A)$.

Assume from now on that this holds: $H^i(X,\A)=0$ for $i \geq
0$. The composition $L^\hdot\pi_X^*R^\hdot\pi^X_*$ also comes
equipped with a canonical adjunction map
$L^\hdot\pi_X^*R^\hdot\pi^X_* \to \Id$; our goal in this Subsection
is to find a way to measure whether it is an isomorphism or not.

To do this, consider the product $X \times_S X$ and equip it with
the coherent algebra $\A^{opp} \boxtimes \A$ ($\A^{opp}$ denotes the
opposite algebra). Since $H^i(X,\A)=0$ for $i \geq 1$, the algebra
$A=H^0(X,\A)$ is flat over $B$, and by the K\"unneth formula we have
$H^0(X \times_S X,\A^{opp} \boxtimes \A) \cong A^{opp} \otimes_B
A$. Denote by $\pi_2:X \times_S X \to X$ the projection onto the
second factor, and let $\pi_{2*}$ be the associated direct image
functor. For any object $\K \in D^-(X \times_S X,\A^{opp} \boxtimes
\A)$, we define the kernel functor $F(\K):D^-(X,\A) \to D^-(X,\A)$
by setting
$$
F(\K)(\F) = R^\hdot\pi_{2*}\left((\F \boxtimes \A) \lotimes_{\A \boxtimes
\A^{opp}} \K\right)
$$
for any $\F \in \D^-(X,\A)$. An obvious example is the identity
functor; we denote the corresponding kernel by $\A_\Delta \in \Shv(X
\times X,\A^{opp} \boxtimes \A)$. This is a sheaf supported on the
diagonal $\Delta \subset X \times_S X$, and for every open subset $U
\subset X$ we have $\A_\Delta(U \times_S U) = \A(U)$ with the
natural $\A(U)$-bimodule structure.

\begin{lemma}
The functor $L^\hdot\pi_X^*R^\hdot\pi^X_*:D^-(X,\A) \to D^-(X,\A)$
is isomorphic to the kernel functor $F(\M)$ with the kernel
$\M^\hdot = \M^\hdot(X/S,\A) = L^\hdot\pi_{X \times_S X}^*A$, where
$A$ is equipped with the natural $A$-bimodule structure.
\end{lemma}

\proof{}(Compare \cite[Appendix D]{Ku}.) For any $A$-bimodule $N$
and any sheaf $\F \in \Shv(X,\A)$, the projection formula gives an
isomorphism
$$
R^\hdot\pi^X_*\F \otimes_A N \cong R^\hdot\pi^X_*\left(\F \lotimes_\A
L^\hdot\pi_X^*N\right),
$$
while adjunction gives a map
$$
L^\hdot\pi_X^*N \to R^\hdot\pi_{1*}L^\hdot\pi^*_{X \times_S X}N,
$$
and the projection formula again gives an isomorphism
$$
\F \otimes_\A R^\hdot\pi_{1*}L^\hdot\pi^*_{X \times_S X}N
\cong R^\hdot\pi_{1*}\left( \F \boxtimes \A \lotimes_{\A \boxtimes
  \A^{opp}} L^\hdot\pi^*_{X \times_S X}N\right).
$$
Composing all three, we obtain a map
$$
R^\hdot\pi^X_*\F \lotimes_A N \to R^\hdot\pi^{X\times_SX}_*\left( \F
\boxtimes \A \lotimes_{\A \boxtimes \A^{opp}} L^\hdot\pi^*_{X
\times_S X}N\right),
$$
which gives by adjunction a base-change map
$$
L^\hdot\pi_X^*\left(R^\hdot\pi^X_*(\F) \lotimes_A N\right) \to
R^\hdot\pi_{2*}\left((\F \boxtimes \A) \lotimes_{\A \boxtimes
\A^{opp}} L^\hdot\pi_{X \times_S X}^*N\right).
$$
We have to prove that it is an isomorphism for $N=A$. More
generally, we will prove that it is an isomorphism for every
$A$-bimodule $N$. Indeed, since every $A$-bimodule has a free
resolution, it suffices to consider the free bimodule $N = A^{opp}
\otimes_B A$. Then the left-hand side is isomorphic to
$H^\hdot(X,\F) \otimes_B \A$, and the right-hand side is isomorphic
to $R^\hdot\pi_{2*}(\F \boxtimes \A)$, which are the same
$\A$-module by the K\"unneth formula.  \endproof

We see that the two functors that we want to compare are kernel
functors, with kernels $\A_\Delta,\M^\hdot(X/S,\A) \in D^-(X \times
X,\A^{opp} \boxtimes \A)$. The adjunction map is given by a map
$\M^\hdot(X,\A) \to \A_\Delta$.

\begin{defn}\label{K}
Denote by $\K^\hdot(X/S,\A) \in D^-(X \times_S X,\A^{opp} \boxtimes
\A)$ the cone of the adjunction map
$$
\M^\hdot(X,\A) = L^\hdot\pi^*A \to \A_\Delta.
$$
\end{defn}

Thus the adjunction map is an isomorphism if and only if
$\K^\hdot(X/S,\A) = 0$; if this happens, $R^\hdot\pi^X_*$ and
$L^\hdot\pi_X^*$ are mutually inverse equivalences of categories. We
will need one refinement of this statement.

\begin{lemma}\label{dimX}
Assume that $H^i(X,\A)=0$ for $i \geq 1$, and that the category
$\Shv(X \times_S X,\A^{opp} \boxtimes \A)$ has homological dimension
$k$. Then if the complex $\K^\hdot(X/S,\A)$ has no homology in
degrees $\geq -(k-1)$, we have $\K^\hdot(X/S,\A)=0$.
\end{lemma}

\proof{} The kernel $\A_\Delta$ is actually a sheaf placed in degree
$0$, while the kernel $\M^\hdot(X/S,\A)$ is trivial in positive
degrees by construction. Therefore under the assumptions of the
Lemma, the connecting map in the exact triangle
$$
\begin{CD}
\K^\hdot(X/S,\A) @>>> \M^\hdot(X/S,\A) @>>> \A_\Delta @>>>
\end{CD}
$$
must be equal to $0$, and we have $\M(X/S,\A) = \A_\Delta \oplus
\K(X/S,\A)$. In other words, the functor
$F=L^\hdot\pi^*R^\hdot\pi_*$ splits into a direct sum of the
identity functor $\Id$ and some other functor $F'$, in such a way
that the natural adjunction map $\ad:F \to \Id$ becomes the projection
$\Id \oplus F' \to \Id$ onto the first summand. Since $H^i(X,\A)=0$
for $i \geq 1$, we have $R^\hdot\pi^X_*L^\hdot\pi_X^* \cong \Id$, so
that the natural map
$$
F = L^\hdot\pi_X^*R^\hdot\pi^X_* \to
L^\hdot\pi_X^*R^\hdot\pi^X_*L^\hdot\pi_X^*R^\hdot\pi^X_* \cong 
F \circ F
$$
is an isomorphism. Moreover, the composition
$$
\begin{CD}
F @>>> F \circ F @>{\ad \circ \id}>> \Id \circ F \cong F
\end{CD}
$$
tautologically coincides with the identity map. Thus
$$
\ad \circ\id:F \circ F \to \Id \circ F \cong F
$$
must be an isomorphism, too. But it vanishes on the direct summand
$$
F' \cong F' \circ \Id \subset (\Id \oplus F') \circ (\Id \oplus F')
= F \circ F.
$$
We conclude that $F' = 0$.
\endproof

\begin{remark}\label{radi}
In all of the above, the scheme structure on $X$ is really
irrelevant; in particular, if we have a finite scheme map $p:X \to
X_1$ which is identical on the level of points -- for example, a
Frobenius map -- then $p_*\K(X/S,\A)$ tautologically coincides with
$\K(X_1/S,p_*\A)$. In fact, one can develop the same theory for
general ringed topological spaces, or for ringed toposes, with some
appropriate coherence conditions. We do not do this since we do not
need this; however, it might be useful in some applications (see
e.g. Remark~\ref{D-aff-rem}). Ringed topological spaces $\langle
X,\A \rangle$ such that $\pi_*$ and $\pi^*$ are equivalences of
categories are sometimes called {\em $\A$-affine}. A well-known
example is the flag variety $X=G/B$ associated to a reductive
algebraic group $G$; in this example, $\A$ is the algebra of
differential operators on $X$ (see \cite{BB}).
\end{remark}

\subsection{Critical lines.}

Assume now, as in Section~\ref{qua}, that the base field $k$ has odd
positive characteristic $\cchar k = p$. For any $a \in \Z/p\Z$,
denote by $S_a \subset S_h=\Spec k[[h,t]]$ the line
$t=ah$. Moreover, consider the polynomial
$$
P(h,t) = t(t^{p-1}-h^{p-1}) = \prod_{a \in \Z/p\Z} (t-ah),
$$
and denote by $S_s = \Spec k[[h,t]]/P(h,t)$ the union of all lines
$S_a \subset S_h$, $a \in \Z/p\Z$. Equivalently, $S_s \subset S_h$
is the fiber of the map $s:S_h \to S$ over the special point $o \in
S = \Spec k[[t]]$. The projection $\sigma:S_h \to \Spec k[[h]]$ is
in particular flat on $S_s \subset S_h$; for any $a \in \Z/p\Z$, it
induces an isomorphism $\sigma:S_a \to \Spec k[[h]]$. The scheme
$S_s$ is a local scheme with special point $o$, and the open
complement $\compl{S}_s = S_s \setminus \{o\}$ decomposes
$$
\compl{S}_s = \coprod_{a \in \Z/p\Z} \compl{S}_a,
$$
where we denote $\compl{S}_a = S_a \setminus \{o\}$. The map
$\sigma$ identifies $\compl{S}_a$ with the point $\Spec k((h))$.

Assume given a symplectic scheme $X/k$ projective over a normal
affine scheme $Y/k$. Fix an ample line bundle $L$ on $X$, and assume
given a twistor deformations $\X$ of the pair $\langle X,L
\rangle$. Moreover, assume that the twistor deformation $\X$ is
exact in the sense of Definition~\ref{tw.ex}, and that
$H^i(\X,\calo_X) = 0$ for $i \geq 1$. Let $\X^\tw_h$, $\Y^\tw_h$ and
$\X^\tw_o$ be as in Subsection~\ref{qua}.

For any scheme $Z/k$ and any complete local $k$-algebra $B$, we will
denote by $Z \whtimes \Spec B$ the completed product,
$$
Z \whtimes \Spec B = \sspec_Z\calo_Z \whotimes_k B,
$$ 
where the product on the right-hand side is completed with respect
to the adic topology on $\calo_Z \otimes_k B$ induced from the adic
topology on $B$. In particular, consider the schemes $\X^\tw_h$,
$\Y^\tw_h$ defined in \eqref{xh} and \eqref{yh} of
Subsection~\ref{qua}. Then we have
$$
\begin{aligned}
\X^\tw_h &= \X^\tw \whtimes_{S^\tw} S_h,\\
\Y^\tw_h &= S_h \whtimes_{S^\tw} \Y^\tw, 
\end{aligned}
$$
where in the first equation, $S_h$ is mapped to $S^\tw$ by the
map $s$. The scheme $\X^\tw_h$ is projective
over $\Y^\tw_h$, and both are flat over $S_h$. Let $\X^\tw_o \subset
\X^\tw_h$ be as in Subsection~\ref{qua},
$$
\X^\tw_o = \sspec_{\X^\tw}\calo_\X^p \otimes_{k[[t^p]]} k[[t]].
$$
Moreover, let $X^\tw_s = \X^\tw_h \times_{S_h} S_s$. We note that in
fact $X^\tw_s \cong X^\tw \whtimes S_s$; in particular, it does not
depend on the twistor deformation $\X$. The flat projection
$\sigma:S_s \to \Spec k[[h]]$ defines a flat projection
$\sigma:X^\tw_s \to X^\tw_h$ of degree $p$, where we denote $X^\tw_h
= X^\tw \whtimes \Spec k[[h]]$. The map $\sigma$ is ramified over
the special point $o \in \Spec k[[h]]$; in fact the fiber $X^\tw_t =
X^\tw_s \times_{\Spec k[[h]]} \{o\}$ -- or equivalently, the
scheme-theoretic intersection $X^\tw_t = X^\tw_s \cap \X^\tw_o$ --
is easily seen to coincide with $X^\tw \times \Spec k[t]/t^p$, and
the map $\sigma$ simply projects this to $X^\tw$. The complement
$\compl{X}^\tw_s = X^\tw_s \setminus X^\tw_t$ decomposes
\begin{equation}\label{spctr}
\compl{X}^\tw_s = \coprod_{a \in \Z/p\Z}\compl{X}^\tw_a,
\end{equation}
where we denote $X^\tw_a = X^\tw_h \times_{S_h} S_a$,
$\compl{X}^\tw_a = \X^\tw_a \setminus \X^\tw$. For any $a \in
\Z/p\Z$, the projection $\sigma$ identifies the scheme $X^\tw_a$
with $X^\tw_h$, while $\compl{X}^\tw_a$ is identified with
$\compl{X}^\tw_h = X^\tw_h \setminus X^\tw_h$. Denote also $Y^\tw_h
= Y^\tw \whtimes \Spec k[[h]]$ and $\compl{Y}^\tw_h = Y^\tw_h
\setminus Y^\tw$; then $X^\tw_h$ -- and therefore, $X^\tw_a$ for any
$a \in \Z/p\Z$ -- is projective over $Y^\tw_h$, and $\compl{X}^\tw_a
\cong \compl{X}^\tw_h$ is projective over $\compl{X}^\tw_h$.

\begin{remark}
We must caution the reader that although the scheme
$\compl{X}^\tw_h$ maps to $X^\tw \otimes_k k((h))$, this map is not
an isomorphism, due to the completion involved in the definition of
the scheme $X^\tw_h$. Explicitly, we have
$$
\compl{X}^\tw_h = \sspec_{X^\tw}\calo_{X^\tw}((h)),
$$
where $\calo_{X^\tw}((h)) = (\calo_{X^\tw} \whotimes_k
k[[h]])(h^{-1})$ is the algebra of Laurent series in $h$ with
coefficients in $\calo_{X^\tw}$. This algebra is much larger than
$\calo_{X^\tw} \otimes_k k((h))$. Since $X$ is only projective over
an affine scheme $Y$, all our construction will really be only
defined over this large algebra. In particular, the map
$\compl{X}^\tw_h \to X^\tw \otimes_k k((h))$ is not even bijective
on the level of points (we note, for example, that every
$k((h))$-valued point $x:\Spec k((h)) \to \compl{X}^\tw_h$ extends
to a map $\Spec k[[h]] \to X^\tw_h$).
\end{remark}

Fix now a sheaf $\calo_h$ on $\X^\tw_h$ provided by
Proposition~\ref{qua.prop}, denote by $\calo_s$ its restriction to
$X^\tw_s \subset \X^\tw_h$, and for any $a \in \Z/p\Z$, denote by
$\calo_a$ its restriction to $X^\tw_a \subset X^\tw_s \subset
\X^\tw_h$. By Proposition~\ref{qua.prop}, $\calo_h$ restricts to
$\calo_{\X}$ on $\X^\tw_o \subset \X^\tw_h$, which in turn restricts
to $\Fr_*\calo_X$ on $X^\tw \subset \X^\tw_o$, where $\Fr:X \to
X^\tw$ is the Frobenius map. Hence $\calo_a$ also restricts to
$\Fr_*\calo_X$ on $X^\tw = X^\tw_a \cap \X^\tw_o \subset \X^\tw_h$.

To help the reader to visualize various schemes that we have
introduced, we note that they fit together into a diagram
\begin{equation}\label{X.diag}
\begin{CD}
\compl{X}^\tw_a @>>> \compl{X}^\tw_s @>>> \compl{\X}^\tw_h\\
@VVV @VVV @VVV\\
X^\tw_a @>>> X^\tw_s @>>> \X^\tw_h\\
@AAA @AAA @AAA\\
X^\tw @>>> X^\tw_t @>>> \X^\tw_o,
\end{CD}
\end{equation}
where all the squares are Cartesian, all the vertical arrows in the
bottom row are closed embeddings, and all the vertical arrows in the
top row are the complementary open embeddings. The whole diagram is
obtained from $\X^\tw_h/S_h$ by base change with respect to the
similar diagram
$$
\begin{CD}
\compl{S}_a @>>> \compl{S}_s @>>> \compl{S}_h\\
@VVV @VVV @VVV\\
S_a @>>> S_s @>>> S_h\\
@AAA @AAA @AAA\\
\{o\} @>>> S_t @>>> S.
\end{CD}
$$
Here $S_h = \Spec k[[t,h]]$ is the (formal) two-dimensional affine
plane, $S \subset S_h$ is the line $h=0$, $S_a \subset S_h$ is the
line $ah=t$, $S_t = \Spec k[t]/t^p$, and $S_s \subset S_h$ is the
union of the lines $S_a$ for all $a \in \Z/p\Z$. We also have a
version of \eqref{X.diag} with $\Y$, $Y$ instead of $\X$, $X$.

All in all, the have the following schemes equipped with algebra
sheaves: $X/k$ with $\calo_X$, $\X^\tw_o/S$ with $\calo_\X$,
$X^\tw_s/S_s$ with $\calo_s$, $X^\tw_t/S_t$ with $\calo_\X/t^p$, and
for any $a \in \Z/p\Z$, $X^\tw_a/S_a$ and
$\compl{X}^\tw_a/\compl{S}_a$ with $\calo_a$. By assumptions on $\X$
and by Lemma~\ref{tlt}, all these algebra sheaves have no higher
cohomology. As in Definition~\ref{K}, denote
\begin{alignat*}{3}
\K^\hdot &= \K(X/k,\calo_X),\qquad
&\K^\hdot_o &= \K(\X^\tw_o/S,\calo_\X),\\
\K^\hdot_t &= \K(X^\tw_t/S_t,\calo_\X/t^p),\qquad
&\K^\hdot_s &= \K(\X^\tw_s/S_s,\calo_s),\\
\K^\hdot_a &= \K(X^\tw_a/S_a,\calo_a),\qquad
&\compl{\K}^\hdot_a &= \K(\compl{X}^\tw_a/\compl{S}_a,\calo_a),
\end{alignat*}
By construction, all these kernels are complexes of coherent sheaves
bounded from above (in fact trivial in positive degrees). We note
that we tautologically have
$$
\K^\hdot(X^\tw,\Fr_*\calo_X) \cong \Fr_*\K^\hdot
$$
(see Remark~\ref{radi}); moreover, since the construction of the
kernel $\K^\hdot(X/S,\A)$ is compatible with flat base change with
respect to $S$, and $\calo_a$ restricts to $\Fr_*\calo_X$ on $X^\tw
\subset X^\tw_a$, the kernel $\K^\hdot_a$ restricts to
$\Fr_*\K^\hdot$ on $X^\tw \times X^\tw \subset X^\tw_a \times_{S_a}
X^\tw_a$. Analogously, $\K^\hdot_s$ and $\K^\hdot_o$ both restrict
to $\K^\hdot_t$.

\begin{defn}\label{crit}
A line $S_a \subset S_h$ is called {\em critical} if the complex
$\compl{\K}^\hdot_a \in \D^-(\compl{X}^\tw_h \times_{k((h))}
\compl{X}^\tw_h)$ is non-trivial. For any point $y \in Y$, a line
$S_a$ is called {\em critical at $y$} if it becomes critical after
replacing $Y$ with the localization $\Spec \calo_{Y,y}$.
\end{defn}

\subsection{Bounds.}
Our goal in this Subsection is to bound the number of critical
lines. We do it in the following way. Firstly, assume that the line
bundle $L$ on $X$ is very ample. Then the line bundle $L \boxtimes
L$ on $X \times X$ is also very ample and defines a closed embedding
$X \times X \to \PP_Y = \Pp^{N_Y} \times Y \times Y$ for some
integer $N_Y > 0$. Using this embedding, we can consider $\K^\hdot$
as a complex of sheaves on $\PP_Y$. Analogously, we have closed
embeddings $X^\tw \times X^\tw \to \PP_{Y^\tw} = \Pp^{N_Y} \times
Y^\tw \times Y^\tw$ and $\compl{X}_h \times_{k((h))} \compl{X}_h \to
\compl{\PP}_h = \Pp^{N_Y} \times \compl{Y}_h \times_{k((h))}
\compl{Y}_h$.

Choose now a closed point $y \in Y^\tw$, let $\PP_y \subset
\PP_{Y^\tw}$ be the fiber over the point $y \times y \in Y \times
Y$, and denote its embedding by $i_y:\PP_y \to \PP_{Y^\tw}$. Denote
also $\compl{\PP}_y = \rho^{-1}(\PP_y) \subset \compl{\PP}_h$, where
$\rho:\compl{Y}_k \to Y^\tw$ is the natural projection; by abuse of
notation, keep the notation $i_y:\compl{\PP}_y \to
\compl{\PP}_h$. For any integer $m$ and any $a \in \Z/p\Z$, denote
\begin{equation}\label{c.yam.eq}
C(y,a,m) = \sum_{\substack{0 \geq l \geq -2\dim X,\\ 0 \leq q
\leq N_Y}} \binom{N_Y}{q} \dim_{k((h))} \Hh^{l+q}\left(
\compl{\PP}_y, i_y^*\compl{\K}^\hdot_a(m-q)\right).
\end{equation}

\begin{lemma}\label{>p}
For any line $S_a$ which is critical at the point $y$, we have
$$
C(y,a,m) \geq p^{\dim X}.
$$
\end{lemma}

\proof{} Replacing $Y$ with an open neighborhood of $y \in Y$ if
necessary, we may assume that $S_a$ is critical everywhere. Then for
some point $x \in \compl{\PP}_y$ with embedding $i_x:x \to
\compl{\PP}_y$, the restriction $i_x^*i_y^*\compl{\K}^\hdot_a$ of
the complex $i_y^*\compl{\K}^\hdot_a$ to the point $x$ is
non-trivial. Moreover, since by construction $\calo_a$ is an Azumaya
algebra on $\compl{X}^\tw_a$, the category of sheaves of
$\calo_a$-modules on $\compl{X}^\tw_a \times_{k((h))}
\compl{X}^\tw_a$ has homological dimension $2\dim X$. Therefore by
Lemma~\ref{dimX} we may assume that this restriction
$i_x^*i_y^*\compl{\K}^\hdot_a$ is non-trivial in some degree $l$
with $0 \geq l \geq -2\dim X$. The same is true for any twist
$\compl{\K}^\hdot_a(m)$. Therefore
\begin{equation}\label{dimx.eq}
\dim_{k((h))}\Hh^l(\compl{\PP}_y,i_{x*}i^*_xi^*_y\compl{\K}^\hdot_a(m))
> 0
\end{equation}
for some $l$, $0 \geq l \geq -2\dim X$, and any integer $m$.
However, the complex $\K^\hdot_a$ is by construction a complex of
sheaves of modules over $\calo_a^{opp} \boxtimes \calo_a$, and on
$\compl{X}^\tw_a \subset X^\tw_a$, the algebra $\calo_a$ is a matrix
algebra: we have $\calo_a~|_{\compl{X}^\tw_a} \cong \End(\E)$ for
some vector bundle $\E$ on $\compl{X}^\tw_a$. Therefore by Morita,
$\compl{\K}^\hdot_a \cong \F^\hdot \otimes (\E^* \boxtimes \E)$
for some complex $\F^\hdot$, and the left-hand side of
\eqref{dimx.eq}, being greater than zero, must be at least $\rk \E^*
\boxtimes \E = p^{\dim X}$. On the other hand, the skyscraper sheaf
$i_{x*}k((h))$ concentrated at the point $x \in \compl{\PP}_y$
admits a Koszul resolution, whose terms are the sheaves
$\calo(-q)^{\oplus\binom{N_Y}{q}}$, $0 \leq q \leq N_Y$. By the
projection formula, $i_{x*}i_x^*i_y^*\compl{\K}^\hdot_a(m) \cong
i_{x*}k((h)) \otimes (i_y^*\compl{\K}^\hdot_a(m))$. Therefore there
exists a spectral sequence which converges to
$\Hh^\hdot(\compl{\PP}_y,i_{x*}i_x^*i_y^*\compl{\K}^\hdot_a)$ and
starts with
$$
\bigoplus_{0 \leq q \leq N_Y}\Hh^\hdot\left(\compl{\PP}_y,
  i_y^*\compl{\K}^\hdot_a(m-q)[q]\right)^{\oplus\binom{N_Y}{q}}.
$$
Comparing this to \eqref{c.yam.eq}, we conclude that
$$
\dim_{k((h))}\Hh^l(\compl{\PP}_y,i_{x*}i^*_xi^*_y\compl{\K}^\hdot_a(m))
\leq C(y,a,m),
$$
which proves the claim.
\endproof

Denote now
\begin{equation}\label{c.ym.eq}
C(y,m) = \sum_{\substack{0 \geq l \geq -2\dim X,\\ 0 \leq q
\leq N_Y}} \binom{N_Y}{q} \dim_k \Hh^{l+q}\left(
\PP_y,i_y^*\Fr_*\K^\hdot(m-q)\right),
\end{equation}
and let $C(y)$ be the number of lines $S_a$, $a \in \Z/p\Z$, which
are critical at $y \in Y$ in the sense of Definition~\ref{crit}.

\begin{prop}\label{est.prop}
Assume that $\K^\hdot_o$ is annihilated by $t^M$ for some integer
$M$. Then for any closed point $y \in Y$ and for any integer $m$, we
have 
\begin{equation}\label{main.est}
C(y)p^{\dim X} \leq MC(y,m)
\end{equation}
\end{prop}

\proof{} By Lemma~\ref{>p}, we have
$$
C(y)p^{\dim X} \leq \sum_{a \in \Z/p\Z}C(y,a,m).
$$
Recall that we have a map $\sigma:\compl{X}_o \to \compl{X}_h$, and
by \eqref{spctr}, the left-hand side is the disjoint union of
$\compl{X}_a$, $a \in \Z/p\Z$.  Comparing \eqref{c.yam.eq} and
\eqref{c.ym.eq}, we see that it suffices to prove that
$$
\dim_{k((h))}\Hh^l\left(\compl{\PP}_y,
i_y^*\sigma_*\compl{\K}^\hdot_s(m)\right) \leq
M\dim_k\Hh^l\left(\PP_y,i_y^*\Fr_*\K^\hdot(m)\right)
$$
for any integers $m$ and $l$.

Now, by definition $i_y^*\sigma_*\K^\hdot_s$ is a complex of
coherent sheaves on the projective space $\Pp^{N_Y}_{k[[h]]}$ over
$\Spec k[[h]]$; $\compl{\PP}_y$ is the fiber of this projective
space over the generic point $\Spec k((h)) \subset \Spec k[[h]]$,
while $\PP_y$ is its fiber over the special point $o \in \Spec
k[[h]]$. Moreover, the complex $i_y^*\sigma_*\K^\hdot_s$ is bounded
from above. Therefore we can apply semicontinuity and deduce that
$$
\dim_{k((h))}\Hh^l\left(\compl{\PP}_y,
i_y^*\sigma_*\compl{\K}^\hdot_s(m)\right) \leq
\dim_k\Hh^l\left(\PP_y,i_y^*\sigma_*i_o^*\K^\hdot_s\right),
$$
where $i_o:\X^\tw_o \times_S \X^\tw_o \subset \X^\tw_h \times_{S_h}
\X^\tw_h$ is the embedding of the fiber over the special point $o
\in \Spec k[[h]]$. By base change, $i_o^*\K^\hdot_s \cong
\K^\hdot_t$, so that the right-hand side is equal to $\dim_k\Hh^l(
\PP_y, i_y^*\sigma_*\K^\hdot_t)$. This can be computed by the
spectral sequence associated to the $t$-adic filtration. This
spectral sequence starts with
$$
\Hh^l(\PP_y,i_y^*\Fr_*\K^\hdot(m)) \otimes_k k[t]/t^p.
$$
By assumption $\K^\hdot_o$ is annihilated by $t^M$; hence so does
$\K^\hdot_t$, being isomorphic to its restriction to $X^\tw_t
\times_{S_t} X^\tw_t \subset \X^\tw_o \times_S \X^\tw_o$. To finish
the proof, it suffices to apply inductively the following obvious
linear-algebraic fact.

\begin{lemma}
For any field $k$, any integer $l$, and any complex $E^\hdot$ of
finitely generated $k[[t]]$-modules, we have $\dim_k H^l(E^\hdot)/t
\leq \dim_k E^l/t$.
\end{lemma}

\proof{} Changing $E^l$ and $E^{l-1}$ if necessary, we may insure
that $H^l(E^\hdot)$ and $\dim_k E^l/t$ stay the same, but the module
$E^l$ is flat over $k[[t]]$. Then
\begin{align*}
\dim_k H^l(E^\hdot)/t &\leq \dim_k \Ker d_l/t = \dim_{k((t))} \Ker
d_l(t^{-1}) \\
&\leq \dim_{k((t))} E^l(t^{-1}) = \dim_k E^l/t,
\end{align*}
where $d_l:E^l \to E^{l+1}$ is the differential in the complex
$E^\hdot$.
\endproof

\section{Proofs.}\label{d.aff}

\subsection{Reduction.}
Assume now that again, the base field $k$ is of characteristic $0$,
and we are given a smooth symplectic scheme $X$ of finite type over
$k$ equipped with a projective birational map $\pi:X \to Y = \Spec
H^0(X,\calo_X)$. Fix a closed point $y \in Y$. Replacing $Y$ and
$\Spec k$ with finite Galois covers, we can assume that the residue
field of the point $y \in Y$ is exactly $k$. Fix an ample line
bundle $L$ on $X$. Choose a projective compactification $\compa{Y}
\supset Y$ of the scheme $Y$, a projective compactification
$\compa{X}$ of the scheme $X$, and a projective birational map
$\pi:\compa{X} \to \compa{Y}$ which extends the given map $\pi:X \to
Y$. Moreover, choose these compactifications in such a way that $L$
extends to an ample line bundle on $\compa{X}$.

By \cite[Corollary 2.8]{K2}, the symplectic form $\Omega_X$ on $X$
is exact in the formal neighborhood of the fiber $\pi^{-1}(y)$ -- or
equivalently, $\Omega_X$ becomes exact after passing to $\wh{X} = X
\times_Y \wh{Y}$, where $\wh{Y}$ is the completion of $Y$ at the
point $y$. In other words, there exist a $1$-form $\wh{\alpha}$ on
$\wh{X}$ such that $\Omega_X = d\wh{\alpha}$. Since
$H^0(\wh{X},\Omega^1_{\wh{X}}) = H^0(\wh{Y},\wh{\pi_*\Omega^1_X})$,
we can choose a $1$-form $\alpha_X \in H^0(X,\Omega^1_X)$ which
approximates $\wh{\alpha}$ to arbitrarily high order at $y \in
Y$. In particular, we can insure that $d\alpha_X$ is non-degenerate
at every point $x \in \pi^{-1}(y) \subset X$. Therefore $d\alpha_X$
is non-degenerate on $\pi^{-1}(U)$ for some Zariski open
neighborhood $U \subset Y$ of the point $y \in Y$. Replacing $Y$
with $U$ and $\Omega_X$ with $d\alpha_X$, we can assume that
$\Omega_X$ is exact.

\begin{remark}
In fact, in \cite[Corollary 2.8]{K2} it is erroneously claimed that
the original $\Omega_X$ itself must be exact over a Zariski-open
neighborhood of the point $y \in Y$. This is not true.
\end{remark}

Fix an exact twistor deformation $\langle \X,\LL \rangle$ associated
to the pair $\langle X,L \rangle$ which is provided by
Lemma~\ref{tw.lemma}. Let $\Y = \Spec H^0(\X,\calo_\X)$. Note that
the canonical bundle $K_\X$ is trivial, so that $H^i(\X,\calo_\X) =
H^i(\X,K_\X) = 0$ for $i \geq 1$ by the Grauert-Riemenschneider
Vanishing Theorem. Moreover, by Lemma~\ref{1-1} $\Y$ is the $t$-adic
completion of a scheme $\wt{Y}$ of finite type over $k[t]$, and
$\wt{Y} \otimes_{k[t]} k$ is an \'etale neighborhood of $y \in Y$.
Therefore, replacing $Y$ with this \'etale neighborhood, we can find
a subalgebra $O \subset k$ of finite type over $\Z$, an affine
variety $\Y_O$ of finite type over $O[t]$, and a variety $\X_O$
which is projective over $\Y_O$ so that the $t$-adic completion of
$\Y_O \otimes_O k$ coincides with $\Y$, we have an isomorphism $\X
\cong \X_O \times_{\Y_O} \Y$, and $Y \cong Y_O \otimes_O k$, where
$Y_O = \Y_O \otimes_{O[t]} O$. Replacing $\Y_O$ with a dense open
subset, we may additionally assume that $\Y_O$ and $\X_O$ are flat
over $O$, $X_O$ is smooth symplectic over $O$, $\Y_O$ is normal, the
image $W_O \subset Y_O$ of the exceptional locus of the map $X_O \to
Y_O$ is also flat over $\Spec O$, and both the $1$-form $\alpha_X$
and the line bundle $\LL$ come from a $1$-form and a line bundle on
$\X_O$. Localizing $O$ and possibly shrinking $\Y_O$ even further,
we may assume as well that the point $y:\Spec k \to \Y$ extends to a
section $y_O:\Spec O \to \Y_O$ of the projection $\Y_O \to \Spec O$,
and that $H^i(\X_O,\calo_{\X_O})=0$ for $i \geq 1$. Finally,
shrinking $\Spec O$ and $\Y_O$ even further, we may assume that the
compactifications $\compa{Y}$, $\compa{X}$ come from schemes
$\compa{Y}_O$, $\compa{X}_O$ which are projective over $O$ and
contain $Y_O$ and $X_O = \X_O \otimes_{\Y_O} Y_O$ as open dense
subsets, while the map $\pi:\compa{X} \to \compa{Y}$ comes from a
map $\pi:\compa{X}_O \to \compa{Y}_O$.

For any maximal ideal $\m \subset O$ with residue field
$k(\m)=O/\m$, we can consider the corresponding special fibers $X_\m
= X_O \times_O k(\m)$, $\X_\m = \X_O \times_O k(\m)$ of $X_O/\Spec
O$, $\X_O/\Spec O$. Then no matter what is the characteristic of the
field $k(\m)$, $X_\m$ is a smooth symplectic variety over $k(\m)$
equipped with a projective birational map $\pi:X_\m \to Y_\m = Y_O
\times_O O/\m$, and $\X_\m$ is an exact twistor deformation of
$X_\m$. Moreover, the section $y_O:\Spec O \to Y_O$ defines a closed
point $y_\m \in Y_\m$.

Consider the scheme $X_O/Y_O$, and let 
$$
\K^\hdot_O = \K^\hdot(X_O,\calo_{X_O}) \in D^b(X_O \times_O X_O)
$$
be as in Definition~\ref{K}. By assumption, the line bundle $L$ on
$\compa{X}_O$ is ample. Replacing it with a multiple, we may assume
that it is very ample. Then $L \boxtimes L$ is very ample on
$\compa{X}_O \times_O \compa{X}_O$, and in particular, it defines a
closed embedding $\compa{X}_O \times_O \compa{X}_O \to \compa{\PP}_Y
= \Pp^{N_Y}_O \times_O \compa{Y}_O \times_O \compa{Y}_O$ for some
integer $N_Y > 0$. Moreover, it also defines a global projective
embedding $\compa{X}_O \times_O \compa{X}_O \to \PP_O = \Pp^N_O$ for
some $N > N_Y$, so that we have a chain of embeddings
$$
\begin{CD}
\compa{X}_O \times_O \compa{X}_O @>>> \compa{\PP}_Y @>>> \PP_O.
\end{CD}
$$
By pushforward, we can treat the complex $\K^\hdot_O$ as a complex
of coherent sheaves on $\PP_Y = \Pp^{N_Y} \times_O Y_O \times_O
Y_O$. Extend it to a complex $\compa{\K^\hdot}_O$ on $\PP_O$ in such
a way that for every $l$, the support of the $l$-th homology sheaf
of the complex $\compa{\K^\hdot}_O$ is the closure of the support of
the $l$-th homology sheaf of the complex $\K^\hdot_O$ (no new
components of the support appear at infinity). Denote by
$\hh^\hdot_O$ the homology sheaves of the complex
$\compa{\K}^\hdot_O$, and denote by $\rk^l_O = \dim (\Supp \hh^l_O)
- \dim\Spec O$ the dimension of the support of the sheaf $\hh^l_O$,
relative over $O$.  Moreover, let $\PP_y \subset \PP_Y \subset
\PP_O$ be the preimage of the $O$-valued point $y_O \times y_O \in
Y_O \times_O Y_O$ under the projection $\PP_Y \to \compa{Y}_O
\times_O \compa{Y}_O$. Then $\PP_y \subset \PP_X$ is a linear
subspace in a projective space; therefore the structure sheaf of
$\PP_y \subset \PP_X$ has a Koszul resolution $\F^\hdot$ by sheaves
of the form $\F^r = \calo(-r)^{\oplus\binom{N-N_Y}{r}}$, $0 \leq r
\leq N - N_Y$.

\bigskip

Now, in order to prove Theorem~\ref{main} below in
Subsection~\ref{pf.sub}, we want to show that for almost all maximal
ideals $\m \subset O$ with residue field $k(\m)$ of positive
characteristic, the canonical quantization of the exact twistor
deformation $\X_\m$ contains at least one line which is not critical
in the sense of Definition~\ref{crit}. To do this, we will use
estimates from Section~\ref{est}. However, we need to re-do them in
a way independent of $\m$ and $\cchar k(\m)$.

\medskip

We first note that shrinking $\Spec O$ if necessary, we may assume
that all the sheaves $\hh^l_O$, $0 \geq l \geq -4\dim X$, are flat
over $O$. Moreover, there exists an integer $m_0$ such that for any
$m_1 > m$, the sheaves $\hh^l_O(m_1)$ have no cohomology, and all
global section modules $H^0(\PP_X,\hh^l_O(m_1))$ are flat over
$O$. Increasing $m_0$, we may assume that for any $m_1 \geq m_0$ we
have
\begin{equation}\label{hlb}
2a_lm_1^{\rk_O^l} \geq \rk_O H^0(\PP_X,\hh^l_O(m_1)),
\end{equation}
where $a_l$ is the leading coefficient in the Hilbert polynomial of
the sheaf $\hh^l_O$. We note that by assumption $\hh^l_O$ is flat
over $O$, so that $\rk^l_O = \dim \Supp \hh^l_O \otimes_O k$;
moreover, we have
$$
\dim \Supp \hh^l_O \otimes_O k = \dim \Supp \hh^l_O \otimes_O
  k\mid_{\PP_Y} \leq \dim X \times_Y X,
$$
since $\K^\hdot(X,\calo_X)$ is by construction supported in $X
\times_Y X$. But by \cite[Lemma 2.11]{K2} the map $X \to Y$ is
semismall -- in other words, $\dim X \times_Y X = \dim X$. We
conclude that $\rk^l_O \leq \dim X$ for any $l$.

\medskip

Secondly, we note that by Lemma~\ref{1-1}, the map $\pi:\X \to \Y$
is one-to-one outside of the special fiber $X \subset \X$. In the
language of kernels, this means that $\K^\hdot(\X,\calo_\X)$ is
supported on the special fiber $X \times X \subset \X \times_S \X)$
-- or, in algebraic language, $\K^\hdot(\X,\calo_\X)$ is annihilated
by $t^M$ for some integer $M \gg 0$ (here, as before, $t$ is the
coordinate on the base $S = \Spec k[[t]]$ of the twistor deformation
$\X/S$). Shrinking $\Spec O$ and $\Y_O$ if necessary, we may assume
that $\K^\hdot(\X_O,\calo_{\X_O})$ is annihilated by $t^M$
everywhere, not just over the generic point $\Spec k \subset \Spec
O$.

\medskip

We are now ready to prove the following result.

\begin{prop}\label{crit.bound}
In the assumptions and notations above, there exists an integer $C$
such that for any maximal ideal $\m \subset O$ with residue field
$k(\m) = O/\m$ of odd positive characteristic $p = \cchar k(\m)$,
the canonical quantization of the twistor deformation $\X_\m$
provided by Proposition~\ref{qua.prop} contains at most $C$ lines
$S_a$, $a \in \Z/p\Z$ which are critical at $y_\m \in Y_\m$ in the
sense of Definition~\ref{crit}.
\end{prop}

\proof{} We will use the estimate of Proposition~\ref{est.prop}; the
reader will easily check that all of its assumptions are
satisfied. The only thing to do is to bound the right-hand side of
\eqref{main.est} -- that is, the integer $C(y,m)$ defined in
\eqref{c.ym.eq} -- in a way independent of $p$. Indeed, use the
resolution $\F^\hdot$ of the structure sheaf of $\PP_y \subset
\PP_O$ to compute the functor $i_y^*$. Then for any $q$, $0 \leq q
\leq N_Y$, the homology sheaves of the complex
$i_y^*(\Fr_*\K^\hdot)(m-q)$ in \eqref{c.ym.eq} can be computed by a
spectral sequence which starts with $(\Fr_*\hh^\hdot_\m) \otimes
\F^\hdot(m-q)$. But we have $\F^r =
\calo(-r)^{\oplus\binom{N-N_Y}{r}}$, and
\begin{align*}
H^\hdot(\compa{\PP}_\m,(\Fr_*\hh^l_\m) \otimes \F^r(m-q))
&= H^\hdot\left(\compa{\PP}_\m, \hh^l_\m \otimes
\Fr^*\F^r(m-q)\right)\\
&= H^\hdot \left(\compa{\PP}_\m, \hh^l_\m(pm -pq
-pr)\right)^{\oplus\binom{N-N_Y}{r}}.
\end{align*}
Now fix $m$ to be any integer such that $m > m_0 +2\dim X +
N_Y$. Then by assumption, $\hh^l_O(p(m-q-r)))$ has no higher
cohomology for any $r \leq N$. Moreover, by construction we have $q
\leq N_Y$ and $r \leq N-N_Y$; therefore by \eqref{hlb}
$$
\dim H^0(\compa{\PP}_\m, \hh^l_\m(pm -pq-pr)) 
\leq 2a_lp^{\rk^l_O}(m-q-r)^{\rk^l_O},
$$
where, since $\hh^l_O$ is flat over $O$, the highest coefficient
$a_l$ of the Hilbert polynomial of the sheaf $\hh^l_\m$ is the
same for all $\m$. Collecting all this together, we see that
\begin{multline*}
\dim H^l(\PP_\m,i_y^*(\Fr_*\K^\hdot)(m-q)) \leq\\
\leq p^{\rk^l_O} \sum_{0 \leq r \leq
  l}2a_{l-r}\binom{N-N_Y}{r}(m-q-r)^{\rk^l_O}.
\end{multline*}
The sum in the right-hand side is a well-defined integer independent
of the maximal ideal $\m \subset O$. It remains to compare this to
\eqref{main.est} and \eqref{c.ym.eq}, and to notice that, since
$\rk^l_O \leq \dim X$, we obtain a bound on the number $C(y)$ of
critical lines independent of $\m$. Summing over all $l$, $0 \geq l
\geq -2\dim X$, and over all $q$, $0 \leq q \leq N_Y$, gives the
desired constant $C$.
\endproof

\begin{remark}
It is perhaps useful to sum up briefly the essential points of the
above proof. We start with a trivial observation: if for some
finitely generated $S_o$-module $E$, the quotient $E/h$ is
annihilated by $t^M$, then, when we consider the decomposition
$$
E(h^{-1}) = \bigoplus_{a \in \Z/p\Z}E_a,
$$
the module $E_a$ is non-trivial for at most $M$ values of $a \in
\Z/p\Z$. To apply this, we need to pass from modules to coherent
sheaves on $X_o \times_{S_o} X_o/S_o$, which requires two main
estimates:
\begin{enumerate}
\item Since $\compl{\K}^\hdot_a$ is a complex of sheaves of modules
  over a matrix algebra of rank $p^{2\dim X}$, every non-trivial
  fiber $i_x^*\compl{\K}^\hdot_a$ of this complex has homology
  spaces of dimension at least $p^{\dim X}$ (this is
  essentially Proposition~\ref{est.prop}).
\item Since the kernel $\K^\hdot$ is supported on the fibered
  product $X \times_Y X \subset X \times X$, and $X \to Y$ is
  semismall, the dimensions of the vector spaces
  $\Hh^\hdot(\compl{\PP}_h,i_y^*\Fr_*\K^\hdot(m))$ are bounded from
  above by $Ap^{\dim X}$, where the constant $A$ does not depend on
  $\m \subset O$ (this is the key part of
  Proposition~\ref{crit.bound}).
\end{enumerate}
The rest of the argument is straightforward homological algebra: we
show that the spaces whose dimension is bounded from below in
\thetag{i} are obtained as an $E^\infty$-term of a spectral
sequence, and the dimensions of the spaces in the $E^0$-term of this
sequence are bounded from above in \thetag{ii}.
\end{remark}

\subsection{Lifting.}\label{pf.sub}
Proposition~\ref{crit.bound} essentially shows that in the
assumptions of Theorem~\ref{main}, and possibly after replacing
$Y$ with an \'etale neighborhood of the point $y \in T$, one can
construct a tilting generator $\E$ on a generic specialization
$X_\m$ of the smooth symplectic variety $X$. To prove
Theorem~\ref{main}, it remains to lift this generator back to
characteristic $0$.

\proof[Proof of Theorem~\ref{main}.]  We are given a smooth
symplectic scheme $X$ of finite type over a field $k$ of
characteristic $0$; we assume that it is equipped with a projective
birational map $\pi:X \to Y = \Spec H^0(X,\calo_X)$. We fix a closed
point $y \in Y$, an ample line bundle $L$ on $X$, and a twistor
deformation $\X$ of the pair $\langle X,L \rangle$ provided by
Lemma~\ref{tw.lemma}.

As explained in the last Subsection, we may choose a subalgebra $O
\subset k$ of finite type over $\Z$, an affine scheme $Y_O$ smooth
and of finite type over $O$, a scheme $X_O$, smooth and symplectic
over $O$ and equipped with a projective birational map $X_O \to
Y_O$, a map $\tau:Y_O \otimes_O k \to Y$ and an $O$-valued point
$y_O:\Spec O \to Y_O$ so that
\begin{enumerate}
\item the map $\tau:Y_O \otimes_O k \to Y$ is \'etale and maps the
  point $y_O \otimes_O k$ to $y \in Y$, and
\item we have $X_O \cong X \times_Y Y_O$, this isomorphism is
  compatible with symplectic forms, and the form on $X_O$ is exact.
\end{enumerate}
Moreover, we may assume that $X_O$ is embedded into a scheme $\X_O$
so that for any maximal ideal $\m \subset O$ with residue field
$k(\m) = O/\m$ of odd positive characteristic, the special fiber
$\X_\m = \X_O \otimes_O k(\m)$ is an exact twistor deformation of
the special fiber $X_\m = X_O \otimes_O O/\m$, and that
$H^i(\X_O,\calo_{\X_O}) = H^i(\X_\m,\calo_{\X_\m})=0$ for $i \geq
1$. By Proposition~\ref{qua}, the twistor deformation $\X_\m$ has a
canonical quantization.  By Proposition~\ref{crit.bound}, we can
assume that for a fixed constant $C$ and any $\m \subset O$ with
$p=\cchar k(\m) > 2$, at most $C$ of the lines $S_a$, $a \in \Z/p\Z$
are critical for the canonical quantization of the deformation
$\X_\m$ in the sense of Definition~\ref{crit}.

Fix an arbitrary maximal ideal $\m \subset O$ such that $p=\cchar
k(\m) > C$, and the field $k(\m)$ is perfect (since $O$ is of finite
type over $\Z$, we can even assume that $k(\m)$ is finite). Then for
some $a \in \Z/p\Z$, the line $S_a$ is not critical for the
quantization of $\X_\m$. This means that the vector bundle $\E$ on
$$
\compl{X}_\m = \sspec(X_\m^\tw,\calo_{X^\tw_\m}((h))),
$$ corresponding to the line $S_a$ in the canonical quantization of
the twistor deformation $\X_\m$ is a tilting generator for the
derived category $D^b_{coh}(\compl{X}^\tw_\m)$ in the sense of
Definition~\ref{tilt}. Since $k(\m)$ is perfect, we can invert its
Frobenius map, so that $X^\tw_\m \cong X_\m$, and $\E$ is actually a
tilting generator of the category $D^b_{coh}(\compl{X}_\m)$. In the
language of kernels, this means that
$\K^\hdot(\compl{X}_\m/k(\m)((h)),\eend(\E)) = 0$.

Denote by $\wh{O}$ the completion of the Laurent series algebra
$O((h))$ with respect to the ideal $\m((h)) \subset O((h))$, denote
by $\wh{Y}_\m$ the completion of the affine variety 
$$
\compl{Y}_O = \sspec(Y_O,\calo_{Y_O}((h)))
$$
along $\compl{Y}_\m \subset \compl{Y}_O$, and denote $\wh{X}_\m =
\wh{Y}_\m \times_{Y_\m} X_\m$. Then $\wh{X}_\m$ is projective over
$\wh{Y}_\m$. Therefore by Grothendieck Algebraization Theorem
\cite[III, Th\'eor\`eme 5.4.5]{EGA}, the vector bundle $\E$ extends from
$\compl{X}_\m$ to the whole $\wh{X}_\m$ if and only if it extends to
the formal scheme neighborhood of $\compl{X}_\m$ in $\wh{X}_\m$. But
since $\E$ is a tilting object -- that is, $\Ext^i(\E,\E)=0$ for $i
\geq 1$ -- the latter is automatic: by standard deformation theory,
$\Ext^2(\E,\E)=0$ implies that $\E$ extends to the infinitesimal
neighborhood $\compl{X}_O \otimes_O O/\m^l$ of arbitrary order $l
\geq 1$, and $\Ext^1(\E,\E)=0$ means that such an extension is
unique. We conclude that $\E$ extends to a vector bundle $\wh{\E}$
on $\wh{X}_\m$. By flat base change, we have
$\Ext^i(\wh{\E},\wh{\E})=0$ for $i \geq 0$.

Now, we know that $\wh{X}_\m$ is projective and of finite type over
$\wh{Y}_\m$. Therefore by \cite[Theorem 1.10]{A} there exists an
\'etale neighborhood $Y_h$ of the point $\compl{y}_\m \in
\compl{Y}_O$ and a scheme $X_h$ projective over $Y_h$ such that
$\wh{X}_\m = X_h \times_{Y_h} \wh{Y}_\m$ and the tilting generator
$\wh{\E}$ comes from a vector bundle $\E_h$ on $X_h$ -- more
precisely, $\wh{E} \cong \tau^*\E_h$, where $\tau:\wh{X}_\m \to X_h$
is the natural projection map. Moreover, we can choose an
$O$-subalgebra $O' \subset \wh{O}$, flat and finite over $O$ and
with fraction field $k'$, so that $Y_h$ and $X_h$ are defined over
$O'((h))$, $\wh{O}$ is faithfully flat over $O'((h))$, the point
$y_\m \subset Y_\m$ extends to an $O'((h))$-valued point $y_h:\Spec
O'((h)) \to Y_h$, and the natural map $Y_h \times_{O'((h))} k'((h))
\to Y \times_k k'((h))$ sends $y_h \times_{O'} k'$ to $y \times_k
k'((h))$ and is \'etale near $y_h \times_{O'((h))} k'$. Replacing
$O'$ and $Y_h$ with \'etale covers of open subsets, we may in
addition assume that the vector bundle $\E_h$ on $X_h$ satisfies
$\Ext^i(\E_h,\E_h)=0$, $i \geq 1$. Therefore the kernel
$\K^\hdot=\K^\hdot(X_h/O'((h)),\eend(\E_h))$ on $X_h
\times_{O'((h))} X_h$ is well-defined, see Definition~\ref{K}. Since
$\K^\hdot \otimes_{O'} \wh{O} =
\K^\hdot(\wh{X}_\m/\wh{O},\eend(\wh{E}))$ reduces to
$\K^\hdot(\compl{X}_\m/k(\m)((h)),\eend(\E)) = 0$ over $\Spec
k(\m)((h)) \subset \Spec \wh{O})((h))$, we have $\K^\hdot
\otimes_{O'} \wh{O} = 0$ by Nakayama Lemma, and since $\wh{O}$ is
faithfully flat over $O'$, this implies that $\K^\hdot=0$. In other
words, $\E$ is a tilting generator for $D^b_{coh}(X_h)$.

Finally, we need to remove the quantization parameter $h$. To do
this, we note that as before, there exists a subalgebra $O_h \subset
k'((h))$ of finite type over $k'$, a scheme $Y'_h$ of finite type
over $O_h$, an $O_h$-valued point $y'_h:\Spec O_h \to Y_h$, and a
scheme $X'_h$ smooth and projective over $Y'_h$ such that $Y'_h
\otimes_{O_h} k'((h)) \cong Y_h$, $y'_h \times_{O_h} k'((h)) = y_h$,
and $X'_h \times_{Y'_h} Y_h \cong X_h$. Moreover, possibly passing
to \'etale covers of open subsets and localizing $O_h$, we may
assume that the natural maps $Y'_h \to \Spec O_h$ and $Y'_h \to Y$
are smooth, and that the tilting generator $\E$ comes from a tilting
generator on $X'_h$. Take a closed point $x \in \Spec O_h$. Its
residue field $K$ is then a finite extension of the field $k'$,
hence of the original field $k$. The fiber $Y_K$ of the scheme
$Y'_h$ over $x \in O_h$ is an \'etale neighborhood of $y \in Y$ and
satisfies the conditions of Theorem~\ref{main}.  \endproof

\begin{remark}
Artin approximation is used two times in our proof of
Theorem~\ref{main}: firstly, in Lemma~\ref{1-1} (to show that
the twistor deformation $\Y$ is algebraic), and secondly, in the
proof itself (to lift the tilting generator from positive
characteristic to characteristic $0$). The second instance is
standard, and there is no way around it. The first one is less
standard. We note, however, that it can in fact be avoided. Indeed,
our argument does not really need the full twistor deformation $\Y$;
at the cost of a slight modification of the notion of a twistor
deformation, one can work just as well with an arbitrary truncation
$\Y/t^N$. After quantization, this only gives the sheaves $\calo_a$
with $a = 0,\dots,N \in \Z/p\Z$; but if $N$ is high enough, this is
sufficient to find a non-critical line.
\end{remark}

\begin{remark}\label{D-aff-rem}
Fedosov quantization was defined originally in characteristic $0$,
not in positive characteristic (see e.g. \cite{BK}). With an
argument similar to ours, one should be able to prove a result
parallel to our Theorem~\ref{main}: for a generic value of
quantization parameter, the manifold $X$ is
$\calo_h(h^{-1})$-affine, where $\calo_h$ is the quantized structure
sheaf. This might even be easier to prove than the corresponding
positive characteristic statement. On the other hand, in examples
such as $X=T^*G/B$, where $G/B$ is the flag variety of a reductive
algebraic group $G$, a stronger statement is true: not only is $X$
$\calo_h(h^{-1})$-affine in the derived category sense, but the
global sections functor is exact, so that even the corresponding
abelian categories are equivalent (see \cite{BB}). It would be
interesting to try to prove a similar statement for general $X$. It
would be also interesting to find exactly the critical values of the
quantization parameter $a$.
\end{remark}

\section{Addenda.}\label{app}

\subsection{$D$-equivalence.}
We will now prove all the additional statements given in
Section~\ref{intr}. We start with the comparison between two
resolutions.

\proof[Proof of Theorem~\ref{K=>D}.]
Assume given two smooth resolutions $X$, $X'$ of a normal
irreducible affine variety $Y$ with trivial canonical bundles $K_X$,
$K_{X'}$. Assume in addition that $X$ admits a closed non-degenerate
$2$-form $\Omega_X$. Fix a closed point $y \in Y$. We note right
away that since $K_X$ is trivial, $Y$ has rational singularities; in
particular, it is a Cohen-Macaulay scheme.

As in the proof of Theorem~\ref{main}, we may replace $Y$ by an
\'etale neighborhood of the point $y$ so that it admits an {\em
exact} non-degenerate form $\Omega_X = d\alpha_X$. Since $X$ and
$X'$ are smooth birational varieties, the form $\alpha_X$ extends to
a $1$-form $\alpha_{X'}$ on $X'$. Its top power
$(d\alpha_{X'})^{\dim X'}$ is a section of the canonical bundle
$K_{X'}$. Since $K_{X'}$ is trivial, the map $X' \to Y$ is
one-to-one over the smooth locus $Y^{sm} \subset Y$, and
$(d\alpha_{X'})^{\dim X'}$ is a global function on $X'$ which is
invertible on $Y^{sm}$. Since $Y$ is normal, $(d\alpha_{X'})^{\dim
X'}$ is non-zero everywhere on $X'$, and $d\alpha_{X'}$ is therefore
a non-degenerate $2$-form on $X'$.

Thus both $X$ and $X'$ are symplectic resolutions of $Y$. By
\cite[Lemma 2.11]{K2}, this means that the maps $\pi:X \to Y$,
$\pi':X' \to Y$ are semismall. In particular, there exists a closed
subvariety $Z \subset Y$ such that $\codim Z \geq 4$, and the fibers
of the maps $\pi$, $\pi^{-1}$ over points in the complement $U = Y
\setminus Z$ are of dimension at most $1$. Since $K_X$ and $K_{X'}$
are trivial, this implies that the graph $\wt{X} \subset X \times_Y
X'$ of the rational map $X \ratto X'$ projects bijectively both to
$X$ and to $X'$ over $U \subset Y$. Thus $X^o=\pi^{-1}(U) \cong
{\pi'}^{-1}(U)$ is a common open subscheme in $X$ and
$X'$. Moreover,
\begin{equation}\label{cm}
H^i(X^o,\calo_{X^o}) = \Hh^i(U,R^\hdot\pi_*\calo_{X^o}) =
H^i(U,\calo_U),
\end{equation}
and since $Y$ is Cohen-Macaulay, this group vanishes for
$i=1,2$.

Let us now trace the steps of the proof of Theorem~\ref{main} for
symplectic schemes $X/Y$, $X'/Y$. First of all, we find a subalgebra
$O \subset k$ of finite type over $\Z$ and $O$-models $Y_O$, $X_O$,
$X'_O$. We obviously can arrange the construction so that the same
$Y_O$ serves for both $X$ and $X'$. Moreover, we can arrange so that
$X^o$ admits a model $X^o_O$, and $H^i(X^o,\calo_{X^o})=0$ for
$i=1,2$ (we can arrange $R^{\geq 1}\pi_*\calo_X=0$, and the model
$Y_O$ is still a Cohen-Macaulay scheme, so that \eqref{cm} applies).
Then we fix an ample line bundle $L$ on $X$, choose a maximal ideal
$\m \subset O$ with residual characteristic $p = \cchar(O/\m)$ and
an element $a \in \Z/p\Z$, and take the deformation $\calo_a$ of the
structure sheaf $\calo_{X_\m}$ which corresponds to a non-critical
parameter $a \in \Z/p\Z$. By Lemma~\ref{tlt}, we have
$\calo_a(h^{-1}) \cong \eend(\E)$ for some tilting vector bundle
$\E$ on $\compl{X_\m}$.

Now, by \cite[Proposition 1.22]{BK3}, there is a family of
one-parameter deformations $\calo_h$ of the structure sheaf
$\calo_X$ which are numbered by elements of the cohomology group
$Q(X)=H^1_{et}(X,\calo_X^*/\calo_X^{*p})$; in \cite{BK3}, they are
called Frobenius-constant quantizations on $X$. In this language,
the deformation $\calo_a$ corresponds to the image of $a[L] \in
\Pic(X)$ under the natural map $\Pic(X) = H^1_{et}(X,\calo_{X}^*)
\to Q(X)$. A similar family of deformations exists for the scheme
$X'_\m$. Since the complement to $X^o_\m \subset X'_\m$ is of
codimension at least $2$, the line bundle $L$ extends to a line
bundle $L'$ on $X_\m'$; denote by $\calo_a'$ the deformation which
corresponds to the image of $a[L'] \in \Pic(X')$ under the map
$\Pic(X') \to Q(X')$. Then the same argument as in Lemma~\ref{tlt}
shows that $\calo_a'(h^{-1}) \cong \eend(\E')$ for some tilting
vector bundle $\E'$ on $\compl{X'_\m}$. On the other hand, since
$H^i(X^o_\m,\calo_{X^o}) = 0$ for $i=1,2$, the deformations
$\calo_a$, $\calo_a'$ coincide on $X^o$. Indeed, while
\cite[Definition 1.21]{BK3} in fact requires control over the
cohomology for $i=1,2,3$, but as the more precise \cite[Proposition
1.18]{BK3} shows, this is only needed to insure the existence of
quantizations; to show that two quantizations with the same
parameter are isomorphic, it is enough to require that the first and
the second cohomology groups are trivial. Thus we obtain vector
bundles $\E$ on $\compl{X}_\m$ and $\E'$ on $\compl{X'}_\m$ such
that $\Ext^i(\E,\E)=\Ext^i(\E',\E')=0$ for $i \geq 1$ and $\E \cong
\E'$ on the common open subset $X^o_\m$. Moreover, since $\eend{\E}$
is isomorphic to $\calo_a(h^{-1})$, we also have
$H^i(X^o_\m,\eend{\E}) = \Ext^i_{X^o_\m}(\E,\E)=0$ for $i=1,2$.

Then we lift $\E$ and $\E'$ to tilting vector bundles on $\wh{X}_\m$
and $\wh{X}'_\m$, and, possibly replacing $Y$ with an \'etale cover,
obtain tilting vector bundles $\E$, $\E'$ on $X$, $X'$. Since
$\Ext^i_{X^o_\m}(\E,\E) = 0$ on $X^o_\m$ for $i=1,2$, the lifting is
unique on $\wh{X^o}_\m$, and the resulting tilting vector bundles
$\E$, $\E'$ are isomorphic on $X^o$. Since the complements $X
\setminus X^o$, $X' \setminus X^o$ have codimension at least $2$
this implies that
$$
R' = \End(\E') = H^0(X',\eend(\E'))=H^0(X^o,\eend(E))=\End(E)=R.
$$
Now, we note that as in the proof of Theorem~\ref{main}, we can
choose $\m \subset O$ and $a \in \Z/p\Z$ so that the resulting
tilting vector bundle $\E$ generates the derived category
$D^b(X)$. Since the line bundle $L'$ is {\em not} ample on $X'$, we
cannot insure the same for $X'$. However, since $D^b(X) \cong
D^b(R\fmod)$, we do know that $R'=R$ has finite global homological
dimension, and we do have a pair of adjoint functors $D^b(R'\fmod)
\to D^b(X)$, $D^b(X) \to D^b(R'\fmod)$. Moreover, their composition
$D^b(R'\fmod) \to D^b(X') \to D^b(R'\fmod)$ is the identity functor,
so that the functor $D^b(R'\fmod) \to D^b(X')$ is a fully faithful
embedding with admissible image. To finish the proof, it suffices to
use the following standard trick.

\begin{lemma}
Assume given an irreducible smooth variety $X$ with trivial
canonical bundle $K_X$ equipped with a birational projective map
$\pi:X \to Y$ to an affine variety $Y$. Then any non-trivial
admissible full triangulated subcategory in $D^b(X)$ coincides with
the whole $D^b(X)$.
\end{lemma}

For the proof we refer the reader, for instance, to \cite[Section
2]{BK2}.
\endproof

\subsection{Positive weights.}
Next, we turn to the situation of a positive-weight $\gm$-action.

\proof[Proof of Theorem~\ref{gm}.] To prove \thetag{i}, we first
note that since the map $\pi:X \to Y$ is semismall by \cite[Lemma
2.11]{K2}, the differential $\xi$ of the $\gm$-action on $Y$ extends
to a vector field on $X$ by \cite[Lemma 5.3]{GK}. Moreover, fix a
relatively very ample line bundle $L$ on $X$; then since
$H^1(X,\calo_X)=0$ by the Grauert-Riemenschneider Vanishing, the
line bundle $L$ admits an action of the vector field $\xi \in
H^0(X,\T(X))$. Replacing $Y$ with its normalization, we may assume
that $Y$ is normal. Then every global section $s \in H^0(X,L^{-1})$
gives an embedding $L \hookrightarrow \calo_X$, which identifies
$H^0(X,L)$ with an ideal $I_s \subset A = H^0(X,\calo_X) =
H^0(Y,\calo_Y)$.  Since $L$ is very ample, the algebra $\bigoplus_i
H^0(X,L^{\otimes i})$ is generated by $H^0(X,L)$, so that we also
have $H^0(X,L^{\otimes n}) = I_s^n \subset A$, and $X$ is isomorphic
to the blow-up of $Y$ in the ideal $I$. Thus to lift the
$\gm$-action on $Y$ to $X$, it suffices to show that $I_s \subset A$
is $\gm$-invariant for an appropriate choice of $s \in
H^0(X,L^{-1})$. Since a subspace in $A$ is $\gm$-invariant if and
only if it is preserved by $\xi$, it suffices to assume that $s$ is
an eigenvector of the vector field $\xi$: $\xi(s)=\lambda s$ for
some $\lambda \in k$.

Indeed, denote by $\m \subset A$ the maximal ideal of the point $y
\in Y$, consider the $A$-module $B=H^0(X,L^{-1})$, and let $\lambda$
be an eigenvalue of the operator $\xi$ on the finite-dimensional
$k$-vector space $B/\m B$. Then since the $\gm$-action on $Y$ is
positive-weight, there exists an integer $N$ such that for any $l
\geq N$, the $\lambda$-eigenspace of $\xi$ on $(\m^l/\m^{l+1})
\otimes (B/\m B)$ is trivial. This means that $\xi - \lambda\id$ is
invertible on $\m^N B \subset B$. Therefore the map $B \to B/\m^N B$
is isomorphic on $\lambda$-eigenspaces, and every
$\lambda$-eigenvector $s_0$ of $\xi$ on the finite-dimensional
$k$-vector space $B/\m^N B$ lifts to a $\lambda$-eigenvector $s \in
B$.

To prove \thetag{ii}, note that by assumption, the tilting vector
bundle $\E$ is defined on $\wh{X} = \wh{Y} \times_Y X$, where
$\wh{Y}$ is the completion of $Y$ at $y \in Y$. Moreover, since $\E$
is tilting, we have $H^1(\wh{X},\eend{\E})=0$, so that $\E$ admits
an action of the vector field $\xi$. Therefore the complete
Noetherian module $M = H^0(\wh{X},\E)$ over the complete Noetherian
algebra $\wh{A} = H^0(\wh{Y},\calo_{\wh{Y}})$ admits an action of
the derivation $\xi:\wh{A} \to \wh{A}$. In order to lift it to a
$\gm$-action, we first have to correct it.

\begin{lemma}
Assume given an Artin $k$-algebra $B$ with ideal $I \subset B$.
Then for every $\gm$-action $\gm \to \Aut(B_0)$ on the quotient
$B_0=B/I$ whose differential lifts to a derivation of $B$, the
composition $\gm \to \gm \to \Aut(B_0)$ with some non-trivial map
$\gm \to \gm$, $\lambda \mapsto \lambda^N$ lifts to a $\gm$-action
on $B$.  Moreover, if $B$ is local with maximal ideal $\m \subset
B$, and $I \subset \m^2$, then we can take $N=1$, and the lifting is
unique.
\end{lemma}

\proof{} (Compare \cite[Section 4]{K2}.) To prove uniqueness, note
that the difference between two liftings is a $\gm$-action on $B$
which is trivial on $\m/\m^2$; therefore its differential $\xi \in
\End_k(B)$ is nilpotent, and being semisimple it must be trivial. To
prove the rest, let $\xi_s$ be the semisimple part in the Jordan
decomposition of the endomorphism $\xi \in \End_k(B')$. The Lie
algebra $\End_k(B)$ acts on the tensor product $B^* \otimes B^*
\otimes B$, a map $a \in \End_k(B)$ is derivation if and only if
$a(m)=0$, where $m \in B^* \otimes_k B^* \otimes B$ is
multiplication $B \otimes_k B \to B$ in $B$. Since the Jordan
decomposition is universal, $\xi(m)=0$ implies $\xi_s(m)=0$, so that
$\xi_s:B \to B$ is also a derivation. Moreover, $\xi_s$ preserves
the line $\Lambda^{\dim I}I \subset \Lambda^{\dim I}(B)$, which
means that $\xi_s$ preserves $I \subset B$. Since $\xi$ is by
assumption semisimple on $B_0$, the difference $\xi-\xi_s$ maps $B$
into $I \subset B$, and $\xi_s=\xi$ on $B_0=B/I$. Denote by $T
\subset \Aut(B)$, $T_0 \subset \Aut(B_0)$ the minimal algebraic
subgroups whose Lie algebras contain $\xi_s$. Then $T_0 \cong \gm$,
$T$ is an algebraic torus, and to finish the proof, it remains to
show that the natural map $T \to T_0$ is surjective. This
immediately follows from the minimality of $T_0 \subset \Aut(B_0)$.
\endproof

Denote by $\m \subset \wh{A}$ the maximal ideal of the local
$k$-algebra $\wh{A}$. Applying this Lemma inductively to the split
square-zero extensions $(A/\m^n A) \oplus (M/\m^n M)$, we see that
the given $\gm$-action on $\wh{A}$ extends to a $\gm$-action on the
module $\wh{M}$. To finish the proof, we apply the following
standard general fact.

\begin{lemma}
Assume given a finitely-generated $k$-algebra $A$ equipped with a
positive-weight $\gm$-action, let $\m \subset A$ be the ideal of
elements of strictly positive weights, and let $\wh{A}$ be the
completion of $A$ with respect to the $\m$-adic topology. Then the
$\m$-adic completion functor is an equivalence between the category
of finitely-generated $\gm$-equivariant $A$-modules and the category
of complete Noetherian $\gm$-equivariant $\wh{A}$-modules.
\end{lemma}

\proof{} Assume given a complete Noetherian $\gm$-equivariant module
$N$. Since $N$ is Noetherian, the quotient $N/\wh{\m}^lN$ is a
finite-dimensional $k$-vector space for every $l$. For
finite-dimensional spaces, a $\gm$-action is the same as a
grading, so that $N/\wh{m}^lN$ splits into a sum of components
$(N/\wh{\m}^lN)^p$ of weight $p$, $p \in \Z$. Consider the space
\begin{equation}\label{fn}
N^f = \bigoplus_p N^p = \bigoplus_p 
\lim_{\gets}\left(N/\wh{\m}^\hdot N\right)^p
\end{equation}
of $\gm$-finite vectors in $N$.  Then $N^f$ is a $\gm$-equivariant
$A$-module, and we claim that $N \mapsto N^f$ is an equivalence
inverse to $M \mapsto \wh{M}$. Indeed, since for any finitely
generated $\gm$-equivariant $A$-module $M$ we obviously have
$\wh{M}/\wh{\m}^l\wh{M} \cong M/\m^l M$, we have $\wh{M}^f \cong
M$. Conversely, for any Noetherian complete $\gm$-equivariant module
$N$ and any $l$, we have a natural map $N^f/\m^lN^f \to
\wh{N}/\m^l\wh{N}$. This map is surjective by definition. Moreover,
since $\m$ has positive weights, for every $p$ the inverse limit in
\eqref{fn} actually stabilizes at some finite level. Therefore every
$a \in N$ is uniquely represented as $a = \lim_{l\to\infty} a^{\leq
l}$, where $a^{\leq l} \in N^f \subset N$ is a finite sum of vectors
$a^i \in N^i$ of weights $i \leq l$. Since $N$ is Noetherian, for
some fixed constant $q \in \Z$ we have $N^l=0$, $l < q$. Then for
every $a \in N$, $a \in \wh{\m}^p N$ means by definition
$$
a = \sum_i m_ia_i
$$
for some $a_i \in N$, $m_i \in \wh{\m}^p$; if $a$ is $\gm$-finite, so
that $a = a^{\leq l}$ for some $l$, this implies
$$
a = \sum_i m_i^{\leq (l-q)}a_i^{\leq l},
$$
so that $a \in \m^pN^f$. Therefore $N^f/\m^p N^f \to N/\wh{\m}^pN$
is bijective.  In particular, the $k$-vector space $N^f/\m N^f$ is
finite-dimensional, so that $N^f$ is finitely-generated over $A$,
and moreover, we have $\wh{N^f} \cong N$.
\endproof

\subsection{Resolution of the diagonal.}\label{pure.sub}
Finally, we turn to cohomological results.

\proof[Proof of Theorem~\ref{diag}.] First of all, by Artin
approximation \cite{A} it suffices to prove the claim under the
additional assumption that $A = H^0(Y,\calo_Y)$ is a {\em complete}
local $k$-algebra. Since $\E$ is a tilting generator for $X$, so it
the dual vector bundle $\E^*$. The product $\E^* \boxtimes \E$ is
then a tilting generator for $X \times X$, so that $D^b(X \times X)
\cong D^b((R^{op} \otimes R)\fmod)$. Under this equivalence, the
diagonal sheaf goes to the tautological $R$-bimodule $R$. Since $X
\times X$ has finite homological dimension, so does the algebra
$R^{op} \times R$. Therefore the bimodule $R$ admits a finite
projective resolution, and to prove the claim, it suffices to shows
that
\begin{enumerate}
\item every indecomposable projective $R$-bimodule is of the form $P'
  \boxtimes P$, where $P'$ is a projective right $R$-module, and
  $P$ is a projective left $R$-module, and
\item every projective left $R$-module $P$ corresponds to a vector
  bundle on $X$ under the equivalence $D^b(X) \cong D^b(R\fmod)$.
\end{enumerate}
The second claim is immediate: every projective module is a direct
summand of $R^N$ for some $N$, and $R^N$ corresponds to $\E^N$. To
prove \thetag{i}, note that since $A$ is complete, so is
$R$. Therefore every indecomposable projective module is a
projective cover of a unique simple $R$-module, the same is true for
$R$-bimodules, and it suffices to prove that every simple
$R$-bimodule is of the form $M' \boxtimes M$, $M'$ a simple right
$R$-module, $M$ a simple left $R$-module. In other words, we may
replace $R$ with the semisimple quotient $R/I$, where $I \subset R$
is the radical of algebra $R$. Then the claim becomes obvious.
\endproof

\proof[Proof of Corollary~\ref{pure}.]
We first prove that the natural map $H^\hdot(X,\Lambda) \to
H^\hdot(F,\Lambda)$ is an isomorphism (here $\Lambda$ is the
coefficient field, such as $\Q_l$ or $\Q$). Indeed, since the
$\gm$-action on $Y$ lifts to an action on $X$ compatible with map
$\pi:X \to Y$, the direct image complex $R^\hdot\pi_*\Lambda$ is
$\gm$-equivariant. By proper base change, it suffices to use the
following.

\begin{lemma}
Assume given a positive-weight $\gm$-action on an affine algebraic
variety $Y$ with fixed point $y \in Y$. Denote by $i_y:y \to Y$ the
embedding. Then for any complex $\F^\hdot$ of constructible sheaves
on $Y$, the natural map $H^\hdot(Y,\F^\hdot) \to i_y^*\F^\hdot$ is
an isomorphism
\end{lemma}

\proof{} We have to show that $H^\hdot(Y,j_!j^*\F^\hdot)=0$, where
$j:U \to Y$ is the embedding of the complement $U = Y \setminus
\{y\}$. By assumption, the algebra $A = H^0(Y,\calo_Y)$ of functions
on $Y$ is positively graded. Let $E = \Proj A^\hdot$, and let
$\wt{Y}$ be the total space of the line bundle $\calo(1)$ on
$E$. Then we have a proper map $\tau:\wt{Y} \to Y$ and an embedding
$\wt{j}:U \to \wt{Y}$. By proper base change, $\tau_*\wt{j}_! \cong
\tau_!\wt{j}_! \cong j_!$, so that it suffices to prove that
$H^\hdot(\wt{Y},\wt{j}_!j^*\F^\hdot)=0$. But since $\F^\hdot$ is by
assumption $\gm$-equivariant, we have $j^*\F^\hdot \cong
\tau^*\G^\hdot$ for some complex $\G^\hdot$ on $E$, and it suffices
to prove that the natural map $H^\hdot(\wt{Y},\tau^*\G^\hdot) \to
H^\hdot(E,\G^\hdot)$ is an isomorphism. By the projection formula,
it suffices to prove that $\tau_*\Lambda_{\wt{Y}}\cong\Lambda_E$.
This claim is local on $E$, so that we may assume $\wt{Y} \cong E
\times \aA^1$, and the claim immediately follows from the K\"unneth
formula.
\endproof

Now, $X$ is smooth; therefore for every integer $l$, all weights in
the cohomology group $H^l(X,\Lambda)$ are $\geq l$. On the other
hand, $F$ is proper; therefore weights in $H^l(F,\Lambda)$ are $\leq
l$. We conclude that for every $l$, $H^l(X,\Lambda)=H^l(F,\Lambda)$
is pure of weight $l$.

Consider now the $\gm$-equivariant cohomology groups
$H^\hdot_{\gm}(X,\Lambda)$ -- that is, the cohomology groups of the
simplicial scheme $X_\idot$, $X_l = X \times \gm^{l}$ with face and
degeneracy maps coming from the $\gm$-action on $X$ (for the theory
of equivariant cohomology, see \cite{BL}). Then we have a spectral
sequence which converges to $H^\hdot_{\gm}(X,\Lambda)$ and starts
with $H^\hdot(X,\Lambda)[u]$, where $u$ is a free generator of
degree $2$. Since $H^\hdot(X,\Lambda)$ is pure, this spectral
sequence collapses, so that $H^\hdot_{\gm}(X,\Lambda) \cong
H^\hdot(X,\Lambda)[u]$. Analogously, for cohomology with compact
support we have $H^\hdot_{\gm,c}(X,\Lambda) \cong
H^\hdot_c(X,\Lambda)[u]$. But since all the $\gm$-fixed points in
$X$ lie within proper subvariety $F \subset X$, the standard
localization theorem (see e.g. \cite{AB}) shows that the natural map
$$
\tau:H^\hdot_{\gm,c}(X,\Lambda) \to H^\hdot_{\gm}(X,\Lambda)
$$
becomes an isomorphism after we invert the parameter $u$. Therefore
it suffices to prove that the image of the map $\tau$ lies in the
$k[u]$-submodule generated by classes of algebraic cycles.  This
immediately follows from the existence of a resolution of the
diagonal provided by Theorem~\ref{diag}. Indeed, by construction all
terms $\E_i$, $\F_i$ in this resolution are direct summands of the
bundle $\E^N$ for some integer $N$. By Theorem~\ref{gm} we may
assume that the tilting generator $\E$ is $\gm$-equivariant. For any
idempotent $P \in \End(\E^N)$, its degree-$0$ component $P^0$ with
respect to the $\gm$-action is also an idempotent, and $\Im P^0 =
\Im P$; therefore we may assume that all the $\E_i$ and $\F_i$ are
$\gm$-equivariant vector bundles. Then we have
$$
\tau(a) = \sum_i (-1)^i \langle [\E_i], a \rangle [\F_i],
$$
where $[\E_i]$, $[\F_i]$ are Chern characters of the bundles $\E_i$,
$\F_i$, and $\langle -,- \rangle$ is the Poincar\'e pairing.
\endproof

\begin{remark}
It is natural to expect that the statement of Corollary~\ref{pure}
holds for $Y$ sufficiently small -- say, local Henselian -- but
without any additional structures such as a group
action. Unfortunately, we have not been able to prove it. In
particular, the above proof does not work: the natural map
$H^\hdot_c(X,\Lambda) \to H^\hdot(X,\Lambda)$ only becomes bijective
after we pass to the equivariant cohomology and localize, and in
ordinary cohomology, it is often $0$ outside of the middle
degree. In practice, a $\gm$-action always exists in the symplectic
case, see \cite[Section 4]{K2}; but in general it is not possible to
prove that it is positive-weight.
\end{remark}

\bigskip

\noindent
{\sc Steklov Math Institute\\
Moscow, USSR}

\bigskip

\noindent
{\em E-mail address\/}: {\tt kaledin@mccme.ru}


{\small
\begin{thebibliography}{BMR}
\bibitem[A]{A} M. Artin, {\em Algebraic approximation of structures
over complete local rings}, Publ. Math. IHES {\bf 36} (1969),
23--58.

\bibitem[AB]{AB} M. Atiyah and R. Bott, {\em The moment map and
equivariant cohomology}, Topology {\bf 23} (1984), 1--28.

\bibitem[Ba]{ba} V. Batyrev, {\em Non-Archimedean integrals and
stringy Euler numbers of log-terminal pairs}, J. Eur.
Math. Soc. {\bf 1} (1999), 5--33.

\bibitem[BB]{BB} A. Beilinson and J. Bernstein, {\em A proof of
Jantzen conjectures}, I. M.  Gel'fand Seminar, 1--50, Adv. Soviet
Math. {\bf 16}, Part 1, AMS, Providence, RI, 1993.

\bibitem[BL]{BL} J. Bernstein and V. Lunts, {\em Equivariant sheaves
and functors}, Lecture Notes in Mathematics, {\bf 1578},
Springer-Verlag, Berlin, 1994.

\bibitem[BK1]{BK} R. Bezrukavnikov and D. Kaledin, {\em Fedosov
quantization in algebraic context}, Moscow Math. J. {\bf 4} (2004),
559-592.

\bibitem[BK2]{BK2} R. Bezrukavnikov and D. Kaledin, {\em McKay
equivalence for symplectic quotient singularities}, Proc. of the
Steklov Inst. of Math., {\bf 246} (2004), 13-33.

\bibitem[BK3]{BK3} R. Bezrukavnikov and D. Kaledin, {\em Fedosov
quantization in positive characteristic}, math.AG/0501247.

\bibitem[BMR]{BMR} R. Bezrukavnikov, I. Mirkovi\'c, and D. Rumynin,
{\em Localization of modules for a semisimple Lie algebra in prime
characteristic}, math.RT/0205144.

\bibitem[BO1]{BO1} A. Bondal and D. Orlov, {\em Semiorthogonal
  decomposition for algebraic varieties}, preprint alg-geom/9506012.

\bibitem[BO2]{BO2} A. Bondal and D. Orlov, {\em Derived categories
  of coherent sheaves}, Proc. ICM 2002 in Beijing, vol. II, Higher
  Ed. Press, Beijing, 2002; 47--56.

\bibitem[BKR]{bkr} T. Bridgeland, A. King, and M. Reid, {\em The
McKay correspondence as an equivalence of derived categories},
J. Amer. Math. Soc. {\bf 14} (2001), 535--554.

\bibitem[Br]{brid} T. Bridgeland, {\em Flops and derived
categories}, Invent. Math. {\bf 147} (2002), 613--632

\bibitem[CP]{CP} C. De Concini and C. Procesi, {\em Quantum groups},
in {\em $D$-modules, representation theory, and quantum groups}
(Venice, 1992), 31--140, Lecture Notes in Math. {\bf 1565},
Springer, Berlin, 1993.

\bibitem[DL]{dl} J. Denef and F. Loeser, {\em Motivic integration,
quotient singularities and the McKay correspondence}, Compositio
Math. {\bf 131} (2002), 267--290.

\bibitem[F]{F} A. Fujiki, {\em On primitively symplectic compact
K\"ahler $V$-manifolds}, in {\em Classification theory of algebraic
and analytic manifolds}, Progress in Math. {\bf 39}, Birkh\"auser,
1983, 71--250.

\bibitem[GK]{GK} V. Ginzburg and D. Kaledin, {\em Poisson
deformations of symplectic quotient singularities}, Adv. Math.  {\bf
186} (2004), 1--57.

\bibitem[Go]{go} I. Gordon and P. Smith, {\em Representations of
symplectic reflection algebras and resolutions of deformations of
symplectic quotient singularities}, Math. Ann. {\bf 330} (2004),
185--200.

\bibitem[EGA]{EGA} A. Grothendieck, {\em \'El\'ements de
G\'eom\'etrie Alg\'ebrique, III}, Publ. Math. IHES {\bf 24}.

\bibitem[H]{H1} D. Huybrechts, {\em Birational symplectic manifolds
and their deformations}, J. Differential Geom. {\bf 45} (1997),
no. 3, 488--513.

\bibitem[K1]{K1} D. Kaledin, {\em On the projective coordinate ring
of a Poisson scheme}, math.AG/0312134, to appear in MRL.

\bibitem[K2]{K2} D. Kaledin, {\em Symplectic singularities from the
Poisson point of view}, math.AG/0310186, to appear in Crelle J.

\bibitem[Ka]{Ka} Y. Kawamata, {\em $D$-equivalence and
$K$-equivalence}, J. Differential Geom.  {\bf 61} (2002), 147--171.

\bibitem[Ku]{Ku} A. Kuznetsov,  {\em Hyperplane sections and derived
  categories},\\ math.AG/0503700.

\bibitem[N]{N1} Y. Namikawa, {\em Deformation theory of singular
symplectic $n$-folds}, Math. Ann. {\bf 319} (2001), 597--623.

\bibitem[R]{R} M. Reid, {\em McKay correspondence}, preprint
alg-geom/9702016.

\bibitem[V1]{vdb1} M. Van den Bergh, {\em Three-dimensional flops
and noncommutative rings}, Duke Math. J. {\bf 122} (2004), 423--455.

\bibitem[V2]{vdb2} M. Van den Bergh, {\em Non-commutative crepant
resolutions}, in {\em The legacy of Niels Henrik Abel}, 749--770,
Springer, Berlin, 2004.

\end{thebibliography}
}
\end{document}